\documentclass[11pt]{amsart}
\usepackage{amsmath,amsthm, amscd, amssymb, amsfonts, mathrsfs,pb-diagram,pb-xy}
\usepackage[all]{xy}
\usepackage{enumitem}
\usepackage{mathrsfs}
\usepackage{srcltx}
\usepackage{hyperref}
\usepackage{multicol}
\usepackage[dvips, dvipsnames, usenames]{color}

\usepackage{ textcomp }

\allowdisplaybreaks

\newcommand{\superqa}[3]{{\bf A}_{#1}(#2|#3)}

\newcommand{\ad}{\operatorname{ad}}
\newcommand{\adc}{\operatorname{ad}_c}
\newcommand{\co}{\operatorname{co}}

\newcommand{\gr}{\operatorname{gr}}

\newcommand{\id}{\operatorname{id}} 
\newcommand{\ord}{\operatorname{ord}}

\newcommand{\GK}{\operatorname{GKdim}}

\newcommand{\ku}{ \Bbbk}
\newcommand{\I}{\mathbb I}
\newcommand{\G}{\mathbb G}
\newcommand{\N}{\mathbb N}
\newcommand{\Z}{\mathbb Z}

\newcommand{\toba}{\mathscr{B}}
\newcommand{\wtoba}{\widetilde{\mathscr{B}}}
\newcommand{\htoba}{\widehat{\toba}}
\newcommand{\bq}{\mathfrak{q}}

\newcommand{\ydh}{{}^{H}_{H}\mathcal{YD}}
\newcommand{\ydkg}{{}^{\ku\Gamma}_{\ku\Gamma}\mathcal{YD}}

\newcommand{\Pc}{\mathcal{P}}

\newcommand{\Zc}{\mathcal{Z}}
\newcommand{\Jc}{\mathcal{J}}

\def\brj{\mathfrak{brj}}

\def\g{\mathfrak{g}}

\def\ufo{\mathfrak{ufo}}

\newcommand{\pre}{\mathfrak{Pre}(V)}
\newcommand{\prefd}{\mathfrak{Pre}_{\operatorname{fGK}}(V)}

\newcommand{\qti}{\widetilde{q}}

\newcommand{\pf}{\begin{proof}}
\newcommand{\epf}{\end{proof}}

\newcommand{\GL}{\operatorname{GL}}

\numberwithin{equation}{section}
\theoremstyle{plain}

\newtheorem{theorem}{Theorem}[section]
\newtheorem{lemma}[theorem]{Lemma}
\newtheorem{definition-theorem}[theorem]{Definition-Theorem}
\newtheorem{prop}[theorem]{Proposition}

\theoremstyle{definition}

\newtheorem{remark}[theorem]{Remark}

\newtheorem{conjecture}[theorem]{Conjecture}

\newtheorem{stepo}{Step}
\newtheorem{stepi}{Step}

\newcommand{\ot}{\otimes}

\begin{document}
\noindent
\title[Finite GK-dimensional pre-Nichols algebras]
{Finite GK-dimensional pre-Nichols algebras of (super)modular and unidentified type}

\author[Angiono]{Iv\'an Angiono}
\address{\noindent FaMAF-CIEM (CONICET) \\ 
Universidad Nacional de C\'ordoba \\
Medina Allende s/n, Ciudad Universitaria \\
5000 C\'ordoba \\
Rep\'ublica Argentina
}
\email{ivan.angiono@unc.edu.ar}
\author[Campagnolo]{Emiliano Campagnolo}
\email{emiliano.campagnolo@mi.unc.edu.ar}
\author[Sanmarco]{Guillermo Sanmarco}
\address{\noindent Department of Mathematics \\ 
Iowa State University \\
Ames, IA 50011 \\
USA
}
\email{sanmarco@iastate.edu}


\thanks{The work of the three authors was partially supported by CONICET and Secyt (UNC)}

\begin{abstract}
We show that every finite GK-dimensional pre-Nichols algebra for braidings of diagonal type with connected diagram of modular, supermodular or unidentified type is a quotient of the distinguished pre-Nichols algebra introduced by the first-named author, up to two exceptions. For both of these exceptional cases, we provide a pre-Nichols algebra that substitutes the distinguished one in the sense that it projects onto all finite GK-dimensional pre-Nichols algebras. We build these two substitutes as non-trivial central extensions with finite GK-dimension of the corresponding distinguished pre-Nichols algebra. We describe these algebras by generators and relations, and provide a basis. 

This work essentially completes the study of eminent pre-Nichols algebras of diagonal type with connected diagram and finite-dimensional Nichols algebra. 
\end{abstract} 

\maketitle

\section{Introduction}

This paper is the third part of a series started in \cite{ASa} and continued in \cite{ACS}, where we contribute to the classification of Hopf algebras with finite Gelfand-Kirillov dimension (abbreviated $\GK$) over an algebraically closed field $\Bbbk$ of characteristic zero. 

Due to the vastness of that problem, further restrictions are usually put in place. For instance, affine Noetherian Hopf algebras with finite $\GK$ have been studied for more than two decades. Also, substantial progress has been made towards the classification of Hopf algebras with small $\GK$. See \cite{BG,BZ,Liu,GZ,WZZ,G-survey} and references therein. 

We focus on \emph{pointed} Hopf algebras with finite $\GK$. Our point of view is inspired by Andruskiewitsch-Schneider program \cite{AS}, which was originally set up for studying finite dimensional pointed Hopf algebras but has proven itself fruitful also on the finite $\GK$ setting, see \cite{AS-crelle, AAH-memoirs,AAH-diag,AAM, And-nilpotent} for example.

Recall that a Hopf algebra $H$ is pointed if the coradical (the sum of all simple subcoalgebras) is just the group algebra of the group-like elements. In this case, the coradical filtration is a Hopf algebra filtration so the associated graded object $\gr H$ is a graded Hopf algebra. If $\Gamma$ denotes the group of group-like elements of $H$,  the Radford-Majid biproduct (or bosonization) \cite{Rad-libro} yields a decomposition $\gr H\simeq R\#\Bbbk\Gamma$, where $R=\oplus_{n\geqslant 0} R^n$ is a coradically graded Hopf algebra in the braided tensor category $\ydkg$ of Yetter-Drinfeld modules over $\Bbbk\Gamma$. At this point we have encountered the three invariants that guide the ongoing problem of classification of pointed Hopf algebras: the group $\Gamma$ of group-like-elements; the Yetter-Drinfeld module $R^1$, called the \emph{infinitesimal braiding}; and the  \emph{diagram} $R$, a Hopf algebra  in the category $\ydkg$. This last invariant serves as the entryway for Nichols algebras into the general problem of classification of pointed Hopf algebras.

Given a Yetter-Drinfeld module $V$ over a group $\Gamma$, we can construct the Nichols algebra $\toba(V)$ of $V$, a coradically graded Hopf algebra in $\ydkg$; via Radford-Majid's biproduct we get a pointed Hopf algebra $\toba(V)\#\Bbbk\Gamma$. This is a distinguished element of the class of all pointed Hopf algebras with coradical $\Bbbk\Gamma$ and infinitesimal braiding $V$ for various reasons. Most importantly for us, if $H$ is in this class, then the subalgebra of the diagram $R$ generated by $V\simeq R^1$ is isomorphic to $\toba(V)$. In other words, diagrams of a pointed Hopf algebra with infinitesimal braiding $V$ are \emph{post-Nichols algebras} of $V$, see \cite{AAR}.

Since we are interested in pointed Hopf algebras with finite $\GK$, a celebrated result of Gromov states that the potential groups of group-like elements, assuming finite generation, necessarily  contain a nilpotent subgroup of finite index. A careful analysis of Nichols algebras over finitely generated nilpotent groups was carried out recently in \cite{And-nilpotent}; in particular, it was explained that much can be afforded with abelian groups. 

In this paper, however, the group $\Gamma$ will only play an implicit role, since most of the theory of Nichols algebras can be expressed in the language of \emph{braided} Hopf algebras. That being the case, it is customary to study families of braided vector spaces rather than Yetter-Drinfeld modules over families of groups. We will focus on braided vector spaces of \emph{diagonal type}, which correspond to semisimple Yetter-Drinfeld modules over finitely generated abelian groups. Such a braided vector space $V$ depends on a square matrix $\bq$ of non-zero elements in the ground field, so we denote $\toba_\bq$ rather than $\toba(V)$.

During the last two decades, remarkable progress has been achieved in the study of Nichols algebras of diagonal type. These Nichols algebras are known to have a restricted PBW basis, which allows to introduce the notion of root systems and certain reflections that give rise to Weyl \emph{groupoids}. Using these tools, Heckenberger classified all braided vector spaces of diagonal type with finite root system. The list contains five families: Cartan type, super type, standard type, (super)modular type and unidentified type \cite{AA17}. Later on, Andruskiewitsch, Angiono and Heckenberger stated the following:

\begin{conjecture}{\cite[1.5]{AAH-memoirs}}\label{conj:AAH}
The root system of a Nichols algebra of diagonal type with finite Gelfand-Kirillov dimension is finite.
\end{conjecture}

While the converse implication holds due to a routinary argument involving PBW basis, the conjecture remains open. However there is strong evidence on its behalf, see \cite{AAH-diag,angiono-garcia}. 
Assuming the validity of this statement, as we will do throughout this paper, we deduce that diagrams of pointed  Hopf algebras with finite $\GK$ and infinitesimal braiding of diagonal type correspond  precisely to finite GK-dimensional post-Nichols algebras of braided vector spaces on Heckenberger's classification \cite{H-classif}. 

As discussed in \cite[\S 2.6]{ASa}, the problem of classifying these post-Nichols algebras can be rephrased in the terms of \emph{pre-Nichols algebras}.
The later language better fits our purposes for two reasons. The first one is rather technical: pre-Nichols algebras of $V$ are certain Hopf quotients of the tensor algebra $T(V)$; this allow us to define pre-Nichols algebras by generators and relations. Thus the family of pre-Nichols algebras becomes partially ordered by projection, determining a poset with minimal element $T(V)$ and maximal one the Nichols algebra. As explained above, we are interested in describing the subposet of pre-Nichols algebras with finite $\GK$; as a initial step we look for minimum elements, called \emph{eminent} pre-Nichols algebras. 

There is a more conceptual reason to study pre-Nichols algebras. For each \emph{finite dimensional} Nichols algebra of diagonal type, the first-named author introduced in \cite{An-dist} the \emph{distinguished} pre-Nichols algebra, which has finite $\GK$ among other crucial features. 
A natural question arises: is the distinguished pre-Nichols algebra eminent? 
For braided vector spaces of  Cartan, super or standard type, that is indeed the case up to very few exceptions, as shown \cite{ASa,ACS}. 
Here we complete the work, studying (super)modular and unidentified types. Summarizing, our main result, after the three papers, is:

\begin{theorem}\label{thm:general}
Let $\bq$ be a braiding matrix such that $\dim\toba_{\bq}<\infty$ and the Dynkin diagram of $\bq$ is connected. Assume that Conjecture \ref{conj:AAH} holds.
\begin{enumerate}[leftmargin=*,label=\rm{(\roman*)}]
\item\label{item:general-1} If $\bq$ is not of type
\begin{itemize}
\item Cartan $A_{\theta}$ or $D_{\theta}$ with $q=-1$,
\item Cartan $A_2$ with $q\in \G_3'$,
\item $\superqa{3}{q}{\{2\}}$ or $\superqa{3}{q}{\{1,2,3\}}$, with $q\in \G_{\infty}$,
\item $\g(2,3)$ with any of the following Dynkin diagram
\begin{align*}
&d_1: \, \xymatrix@C=30pt{\overset{-1}{\circ} \ar@{-}[r]^{\xi^2} &\overset{\xi}{\circ}\ar@{-}[r]^{\xi}\ &\overset{-1}{\circ}}, &
&d_2: \, \xymatrix@C=30pt{\overset{-1}{\circ} \ar@{-}[r]^{\xi} &\overset{-1}{\circ}\ar@{-}[r]^{\xi}\ &\overset{-1}{\circ}},
\end{align*}
\end{itemize}
then the distinguished pre-Nichols algebra $\wtoba_{\bq}$ is eminent.
\item\label{item:general-2} If $\bq$ is of type $\superqa{3}{q}{\{2\}}$ or $\superqa{3}{q}{\{1,2,3\}}$, there is an eminent pre-Nichols algebra $\htoba_{\bq}$ which is a braided central extension of $\wtoba_{\bq}$ by a polynomial algebra in one variable.

\item\label{item:general-3} If $\bq$ is of type Cartan $A_2$ with $q\in \G_3'$ or $\g(2,3)$ with diagram $d_1$ or $d_2$, then $\htoba_{\bq}$ is a braided central extension of $\wtoba_{\bq}$ by a $q$-polynomial algebra in two variables.

\end{enumerate}
\end{theorem}

For braidings of Cartan type $A_{\theta}$ or $D_{\theta}$ with $q=-1$, there are candidates for eminent pre-Nichols algebras. The $\GK$ of these candidates is not known yet, see \cite[\S 5]{ASa}.

\pf
For \ref{item:general-1}, the Cartan case follows by \cite[Theorem 1.3 (a)]{ASa}; super and standard types, by \cite[Theorem 1.2 (1)]{ACS}; the remaining cases are covered by Theorem \ref{thm:distinguished-is-eminent}.

Now \ref{item:general-2} follows from \cite[Theorem 1.2 (2), (3)]{ACS}. Finally Propositions \ref{prop:g(2,3)-case-a} and \ref{prop:g(2,3)-case-b} 
together with \cite[Theorem 1.3 (b)]{ASa} prove \ref{item:general-3}.
\epf

\medbreak

The structure of the paper is the following. In Section  \ref{sec:preliminaries} we fix the notation and recall general aspects of the theory of Nichols algebras for later use. Section \ref{sec:def-relations} is devoted to complete the proof of Theorem \ref{thm:general} \ref{item:general-1}, with a case-by-case analysis of the behaviour of the generators for the defining ideals of Nichols algebras of diagonal type, when considered as elements of pre-Nichols algebras with finite $\GK$. Finally, in Section \ref{sec:g23-exceptions} we introduce and study the eminent pre-Nichols algebras for the two exceptional cases of type $\g(2,3)$ described in Theorem \ref{thm:general} \ref{item:general-3}.

\section{Preliminaries}\label{sec:preliminaries}
In this section we recall definitions and general aspects of the theory Nichols algebras, with particular interest on the diagonal case. We refer to \cite{Tk1} for the definition of braided vector spaces and braided Hopf algebras, and to \cite{Rad-libro,A-leyva,AA17,HS} for unexplained  terminology regarding Hopf and Nichols algebras. For preliminaries on $\GK$ see \cite{KL}.

\subsection*{Notations}
For each $\theta\in\N$ we set $\I_{\theta}=\{1,\dots ,\theta\}$;  we shall write $\I$ when no confusion is possible. 
The canonical basis of the free abelian group $\Z^{\I}$ is denoted $\{\alpha_i: i\in\I\}$; an element $\sum_{i\in\I} b_i\alpha_i$ will be denoted as
$1^{b_1}2^{b_2}\cdots \theta^{b_{\theta}}$.
Given $\beta=1^{b_1}2^{b_2}\cdots \theta^{b_{\theta}}$ and $\gamma=1^{c_1}2^{c_2}\cdots \theta^{c_{\theta}}$, we say that $\beta\le \gamma$ if and only if $b_i\le c_i$ for all $i\in\I$; this defines a partial order on $\Z^{\I}$ that will be used without further mention.

We work over an algebraically closed field $\ku$ of characteristic $0$. The subgroup of $\ku^{\times}$ consisting of $N$-th roots of unity is denoted by $\G_N$, and $\G_ {N}'$ denotes the subset of those of order $N$. The set of all roots of unity is $\G_{\infty}$. 

If $A$ is a $\Z$ graded algebra, we denote by $A^n$ the  degree $n$-component.  The subspace of primitive elements of a (braided) Hopf algebra $H$ is $\Pc(H)$. 

\subsection{Braided action and braided commutator}
Any braided Hopf algebra $\toba$ admits a \emph{left adjoint representation} $\adc: \toba\to\operatorname{End} \toba$,
\begin{align*}
(\adc x)y&=m(m \ot S)(\id\ot c)(\Delta \ot \id)(x \ot y), && x,y\in \toba.
\end{align*}
Also, the \emph{braided bracket} $[\cdot, \cdot]_c:\toba \ot \toba \to \toba$ is the map given by 
\begin{align*}
[x,y]_c &=m(\id-c) (x\ot y), && x,y\in \toba.
\end{align*}

\subsection{Eminent pre-Nichols algebras}
In a nutshell, a \emph{pre-Nichols algebra} of a braided vector space $(V,c)$ is a $\N_0$-graded braided connected Hopf algebra that contains $(V,c)$ as the degree-$1$ component and is generated (as an algebra) by $V$. One can show that there exists a unique pre-Nichols algebra that is moreover \emph{strictly graded}, which means that the subspace of primitive elements coincides with $V$; this is the \emph{Nichols algebra} of $(V,c)$.

There is a more concrete interpretation of these definitions in which one uses the braiding of $V$ to induce in the tensor algebra $T(V)=\oplus_{n\geqslant 0}V^{\ot n}$ a structure of braided Hopf algebra  such that the elements of $V$ are primitive.
Then one can show that there exists a (unique) maximal homogeneous Hopf ideal $\Jc(V)$ of $T(V)$ generated by elements of degree $\ge 2$, thus the Hopf quotient $\toba(V):=T(V)/\Jc(V)$ is the 
Nichols algebra of $V$. 

In this concrete interpretation, a pre-Nichols algebra of $V$ is a Hopf quotient $\toba$ of $T(V)$ by a braided homogeneous Hopf ideal contained in $\Jc(V)$. Hence, the identity of $V$ extends to canonical epimorphisms
$T(V) \twoheadrightarrow \toba \twoheadrightarrow \toba(V)$. 
The family $\pre$ of pre-Nichols algebras of $V$ is naturally equipped with a partial order. Namely,  $\toba_1 \leq \toba_2$ if $\ker (T(V) \twoheadrightarrow \toba_1) \subseteq \ker (T(V) \twoheadrightarrow \toba_2)$. Thus $T(V)$ is minimal and $\toba(V)$ is maximal.

\medbreak
Assume now that $V$ is such that $\GK \toba(V)< \infty$. Inside $\pre$ we get the (non-empty) subposet $\prefd$ containing all finite GK-dimensional pre-Nichols algebras. Understanding this subposet is a crucial problem in the classification of pointed Hopf-algebras with finite GK-dimension, see \cite{ASa} for details. As a first step towards this direction, in \emph{loc.~cit.} the authors introduced the notion of \emph{eminent} pre-Nichols algebra of $V$, which is a minimum of $\prefd$. 
The existence of such minimal objects is not warrantied, see \cite{ASa} for concrete examples related to Lie algebras. However, for the family of braidings of diagonal type with connected diagram, this problem have been addressed in \cite{ASa, ACS}. We will get back to this in \S \ref{sec:def-relations}.

\subsection{Nichols algebras of diagonal type}\label{subsec:diagonal}

A matrix $\bq=(q_{ij})_{i,j \in \I}$ with entries in $\Bbbk^\times$ is called a \emph{braiding matrix}, since it gives raise to a braided vector space $(V,c^\bq)$: fixed a basis $(x_i)_{i \in\I}$ of $V$, $c^{\bq}\in \GL(V\ot V)$ is determined by:
\begin{align}\label{eq:diag-type-def}
c^{\bq}(x_i\ot x_j)&=q_{ij}x_j\ot x_i, & i,j&\in\I.
\end{align}
Such a braiding is called \emph{of diagonal type}. The (Dynkin) diagram of $\bq$ is a decorated graph. The set of vertices is $\I$, each vertex labelled with $q_{ii}$. There is an edge between vertices $i\ne j$ if and only if $\qti_{ij}:=q_{ij}q_{ji}\ne 1$, in such case the edge is labelled with this scalar.

The braiding matrix $\bq$ induces a $\Z$-bilinear form $\bq \colon \Z^{\I}\times\Z^{\I}\to \Bbbk^\times$
defined on the canonical basis by $\bq(\alpha_j,\alpha_k) :=q_{jk}$, $j,k \in\I$.
For $\alpha,\beta  \in \Z^{\I}$ and $i\in \I$,  we set
\begin{align}\label{eq:notation-qab}
q_{\alpha\beta} &= \bq(\alpha,\beta), & q_{\alpha} &= \bq(\alpha,\alpha),&
N_{\alpha} &= \ord q_{\alpha}, & N_{i} &= \ord q_{\alpha_i}.
\end{align}

In the diagonal setting, the Nichols algebra has a compatible $\N_0$-grading. Indeed, the tensor algebra $T(V)$ itself becomes a $\N_0^{\I}$-graded algebra by the rule $\deg x_i=\alpha_i$, $i\in\I$; moreover, the braided structure induced from that of $V$ is homogeneous:
\begin{align*}
c(u\ot v) &= q_{\alpha\beta} \, v\ot u, && u\in T(V)^{\alpha}, \, v\in T(V)^{\beta}, \, \alpha,\beta\in\N_0^{\I}.
\end{align*}
Since the braided coalgebra structure only depends on the braiding, $T(V)$ becomes a $\N_0^{\I}$-graded braided Hopf algebra. Furthermore, the defining ideal $\Jc_{\bq}=\Jc(V)$  turns out to be homogeneuos, thus the Nichols algebra $\toba_{\bq}=\toba(V)$ is a $\N_0^{\I}$-graded braided Hopf algebra.

\medbreak

If $\toba$ is a $\N_0^{\I}$-graded pre-Nichols algebra of $\bq$, the following equalities hold:
\begin{align}
\label{eq:braided-commutator-right-mult}
[u,vw]_c &= [u,v]_c w + q_{\alpha \beta} \, v [u,w]_c,
\\
\label{eq:braided-commutator-left-mult} 
[uv,w]_c &= q_{ \beta \gamma} \, [u,w]_c v + u[v,w]_c,
\\
\label{eq:braided-commutator-iteration}
\big[ [u,v]_c, w \big]_c &= \big[ u, [v, w]_c \big]_c 
- q_{\alpha \beta} \, v [u,w]_c 
+ q_{ \beta \gamma} \, [u,w]_c v,
\end{align}
for all homogeneous elements $u\in \toba^{\alpha}$,  $v \in \toba^{\beta}$, $w \in \toba^{\gamma}$. 

Given at least two indexes $i_1,\dots, i_k\in\I$, we denote 
\begin{align}\label{eq:iterated-adjoint}
x_{i_1\cdots i_k}:=(\adc x_{i_1})x_{i_2\cdots i_k}=x_{i_1}x_{i_2\cdots i_k}-q_{i_1i_2}\cdots q_{i_1i_k}x_{i_2\cdots i_k}x_{i_1} 
\end{align}
as an element in the tensor algebra or any pre-Nichols algebra.

\subsubsection{Classification of arithmetic braidings}\label{subsubsec:classification}
In this work we only consider \emph{arithmetic} braiding matrices $\bq$, which are those with connected diagram and finite generalized root system $\Delta^{\bq}$; that is, those Nichols algebras admitting a (restricted) PBW basis with finite set of generators, so $\Delta^{\bq}_+$ is the set of degrees of a set of generators. This is precisely the class that was classified in \cite{H-classif}, and includes all connected braiding matrices with finite dimensional Nichols algebra. 
Next we recall two results regarding the shape of the diagrams for arithmetic braiding matrices $\bq=(q_{ij})_{i,j\in\I_{\theta}}$.

\begin{lemma}\label{lem:rank3}\cite[Lemma 9 (ii)]{H-classif} 
If $\theta=3$ and the root system is finite, then $\qti_{12}\qti_{13}\qti_{23}=1$ and $(q_{11}+1)(q_{22}+1)(q_{33}+1)=0$. 
Moreover, if $q_{22}, q_{33} \ne -1$ then $q_{22}\qti_{12}=q_{33}\qti_{13}=1$.\qed
\end{lemma}

\begin{lemma}\label{lem:cycles} \cite[Lemma 23]{H-classif}
Assume that $\bq$ has finite root system. Then the diagram of $\bq$ does not contain cycles of length larger than $3$.\qed
\end{lemma}

The classification provided by \cite{H-classif} consists on several tables containing the Dynkin diagram of all arithmetic braidings. Later on, an organization from a Lie-theoretic perspective was achieved in \cite{AA17}, giving rise to five families: Cartan, super, standard, (super)modular, and unidentified. The task of finding eminent pre-Nichols for the first three families was achieved in \cite{ASa, ACS}, up to two exceptions. Here we focus on the remaining families.

\subsubsection{Defining relations}\label{subsubsec:presentation}
An explicit presentation of the Nichols algebras for arithmetic braidings was achieved in \cite[Theorem 3.1]{An-crelle}, which serves as an implicit guidance throughout this work and thus deserves a brief review. That result consists in a list of $29$ homogeneous relations, each of them accompanied by a very specific set of conditions on the entries of the braiding matrix $\bq$ that determine whether or not the relation needs to be included in the presentation of $\toba_\bq$. 

A conceptual analysis of these relations (and their genesis) yields a separation in three categories: quantum Serre relations,
generalizations of these in up to four generators $x_i$, and the so-called \emph{powers of root vectors}. 
Throughout this work, we will refer to this particular set of generators for the ideal $\Jc_\bq$ simply as \emph{the} presentation of $\toba_\bq$.

\subsubsection{Finite GK-dimensional Nichols algebras}\label{subsubsec:Nichols-finite-GK}
Notice that some arithmetic braidings depend on parameters that can take any but a small number of non-zero values in the ground field. In all such cases, the corresponding Nichols algebra is finite dimensional precisely when these parameters are roots of unity; otherwise, it is just finite GK-dimensional. This work and the prequels  \cite{ASa, ACS} focus on braiding matrices with finite dimensional Nichols algebras. The task of finding eminent pre-Nichols algebras when these parameters are not roots of unity will be treated in a sequel.

From this perspective, Conjecture \ref{conj:AAH} states that the classification of Nichols algebras of diagonal type with connected diagram and finite $\GK$ is precisely the one in \cite{H-classif}. Several steps towards proving the conjecture have been achieved. Namely, it is known to be true for braidings of rank $\theta=2, 3$ and of Cartan type, \cite{AAH-memoirs, angiono-garcia}
Some proofs in this work and the prequels \cite{ASa, ACS} assume that the Conjecture holds. However, the  majority of those proofs belong to the realm where the Conjecture is known to be true.

\smallbreak 

As we assume $\GK\toba_{\bq}<\infty$, by \cite[Lemma 20]{Rosso} for each $i\ne j\in\I$ there exists $n\in \N_0$ such that $(\ad_c x_i)^{n+1}x_j=0$.
Then we set
\begin{align*}
	c_{ij}^{\bq}:= -\min \{n\in\N_0: (\ad_cx_i)^{n+1} x_j = 0\}=-\min \{n\in\N_0: (n+1)_{q_{ii}}(1-q_{ii}^n\qti_{ij})=0\}.
\end{align*}
Set also $c_{ii}^{\bq}=2$. Then $C^{\bq}:=(c_{ij}^{\bq})$ is a generalized Cartan matrix and one of the key ingredients in the definition of the Weyl groupoid of $\bq$, cf. \cite{H-Weyl gpd,HY-groupoid}.

\smallbreak 
We end this subsection with two results that will be use several times in \S \ref{sec:def-relations}.

\begin{lemma}\label{lem:1connected} \cite[Proposition 4.16]{AAH-memoirs}
If $W$ is a braided vector space of diagonal type with diagram $\xymatrix@C=30pt{\overset{1}{\circ} \ar@{-}[r]^{p}\ &\overset{q}{\circ}}$, $p\ne 1$, then $\GK \toba(W)=\infty$.\qed
\end{lemma}

\begin{lemma}\label{lem:subspace-primitives} \cite[Lemma 2.8]{ASa}
Let $\toba$ be a graded braided Hopf algebra. If $W$ is a braided vector subspace of $\Pc(R)$, then $\GK \toba(W) \le \GK R$.\qed
\end{lemma}

\subsection{Pre-Nichols algebras of diagonal type}\label{subsec:pre-Nichols-diagonal}
We collect  preliminaries and notation regarding pre-Nichols algebras for later use.

\subsubsection{Distinguished pre-Nichols algebras}\label{subsubsec:distinguished}
Given a connected $\bq$ with finite dimensional $\toba_\bq$, the \emph{distinguished pre-Nichols algebra} $\wtoba_{\bq}$ of $\bq$ is the quotient  of the tensor algebra by the ideal that results from $\Jc_\bq$ by removing certain powers of root vectors and including quantum Serre relations that could formerly be deduced from the power of root vectors. 
This is a key tool in our search for eminent pre-Nichols algebras, since it is designed to admit a PBW basis with the same set of generators as that of $\toba_{\bq}$, and satisfies $\GK\wtoba_{\bq}<\infty$, see \cite{An-dist}.
The powers of root vectors removed from the presentation correspond to \emph{Cartan roots}, the subset of $\Delta^{\bq}_+$ of those PBW generators with infinite height in $\wtoba_{\bq}$.

\subsection{Extensions and Hilbert series of graded braided Hopf algebras}\label{subsec:extensions}
We briefly introduce a tool that will be crucial in \S \ref{sec:g23-exceptions}. For more details, see \cite[\S 2.6]{ACS}.
Following \cite[\S 2.5]{AN}, a sequence of morphisms of braided Hopf algebras 
$\ku \rightarrow A \overset{\iota}{\to} C \overset{\pi}{\to} B \rightarrow \ku $
is an \emph{extension of braided Hopf algebras}  if $\iota$ is injective, $\pi$ is surjective, $\ker \pi = C\iota(A^+)$ and $A=C^{\,\co\pi}$.
In this case, we just write $A \overset{\iota}{\hookrightarrow} C \overset{\pi}{\twoheadrightarrow} B$.

In our examples, $C$ will be connected (i.e. the coradical of $C$ is $\Bbbk$).  By \cite[3.6]{A+}, to get extensions of $C$ it is enough to consider  a surjective  braided Hopf algebra morphism $C \overset{\pi}{\twoheadrightarrow} B$ and set $A=C^{\,\co\pi}$. This construction will be enough for our purposes.

\subsubsection{Hilbert series}
The Hilbert series of  a $\N_0^{\theta}$-graded object $U$ with finite-dimensional homogeneous components is
\begin{align*}
\mathcal{H}_{U} = \sum_{\alpha\in\N_0^{\theta}} \dim U_{\alpha} \, t^{\alpha}
\in \N_0[[t_1,\dots,t_{\theta}]],
\end{align*}
where $t^{\alpha}=t_1^{a_1}\cdots t_{\theta}^{a_{\theta}}$ for $\alpha=(a_1,\cdots, a_{\theta})$.
If $U'$ is $\N_0^{\theta}$-graded object, we say that $\mathcal{H}_U\le \mathcal{H}_{U'}$ if 
$\dim U_{\alpha} \le \dim U'_{\alpha}$ for all $\alpha\in\N_0^{\theta}$.

The main reason for introducing these concepts is the following result.

\begin{lemma}\cite[Lemma 2.4]{ACS}
Fix a Hopf algebra $H$ with bijective antipode. If $A \overset{\iota}{\hookrightarrow} C \overset{\pi}{\twoheadrightarrow} B$ is a degree-preserving extension of $\N_0^{\theta}$-graded connected Hopf algebras in $\ydh$ with finite-dimensional homogeneous components, then $\mathcal{H}_C= \mathcal{H}_A\mathcal{H}_B$.
\end{lemma}

\section{Defining relations and finite GK-dimensional pre-Nichols algebras}\label{sec:def-relations}
Throughout this section we assume that $\bq=(q_{ij})_{i,j\in\I_\theta}$ is a braiding matrix with connected Dynkin diagram such that $\dim\toba_{\bq}<\infty$. In particular the root system is finite (see \cite[\S 3]{H-Weyl gpd}) and each $q_{ii}$ is a root of unity, say of order $N_i$ (necessarily $N_i\geq 2$ by Lemma \ref{lem:1connected}). Let $V^\bq$ be the braided vector space  with basis $(x_i)_{i \in \I_\theta}$ and braiding $c^\bq(x_i\ot x_j)=q_{ij} x_j \ot x_i$.

\begin{theorem}\label{thm:distinguished-is-eminent}
If $\bq$ is not of type
\begin{itemize}
\item Cartan $A_{\theta}$ or $D_{\theta}$ with $q=-1$,
\item $A_2$ with $q\in \G_3'$,
\item $\superqa{3}{q}{\{2\}}$ or $\superqa{3}{q}{\{1,2,3\}}$, with $q\in \G_{\infty}$,
\item $\g(2,3)$ with any of the following Dynkin diagram
\begin{align*}
d_1: \, \xymatrix@C=30pt{\overset{-1}{\circ} \ar@{-}[r]^{\xi} &\overset{-1}{\circ}\ar@{-}[r]^{\xi}\ &\overset{-1}{\circ}},&
&d_2: \, \xymatrix@C=30pt{\overset{-1}{\circ} \ar@{-}[r]^{\xi^2} &\overset{\xi}{\circ}\ar@{-}[r]^{\xi}\ &\overset{-1}{\circ}},
\end{align*}
\end{itemize}
then the distinguished pre-Nichols algebra $\wtoba_{\bq}$ is eminent.
\end{theorem}
\pf
If $\bq$ is either of Cartan, super or standard type, then the proof follows by \cite[Theorem 1.3]{ASa} together with \cite[Theorem 1.2]{ACS}.
Hence we reduce to the cases in which $\bq$ is of types either modular, supermodular or unidentified.

The presentation of the Nichols algebra of $\bq$ given in \cite[Theorem 3.1]{An-crelle} consists on a list of $29$ relations, each of them accompanied by specific conditions on the entries of $\bq$ that determine whether or not the relations needs to be included. Following the procedure in \S\ref{subsubsec:distinguished}, we get a set of relations that give a presentation of the distinguished pre-Nichols algebra $\wtoba_{\bq}$. In the prequels \cite{ASa, ACS} we determined sufficient conditions on $\bq$ to assure that some of these relations hold in any finite $\GK$ pre-Nichols algebra of $\bq$, under some mild assumptions. For these relations the only remaining task is to ensure their validity without any assumption, which is achieved in Lemma 3.1. Finally, in \S \ref{subsec:new-relations} we deal with the relations that where not considered in the prequels.
\epf

\subsection{Relations already considered}\label{subsec:rels-old}
In this subsection we extend some results established in \cite{ASa,ACS} dropping superfluous assumptions. The organization goes as follows. Each relation is studied in a different item of Lemma \ref{lem:previous-lemmas}, where we first fix the elements of $\I_\theta$ that support the relation and then write down the conditions on $\bq$ that \cite[Theorem 3.1]{An-crelle} requires for including this relation in the presentation of $\toba_\bq$. If further hypothesis on $\bq$ are needed, they are included in a different sentence. All such relations are $\N_0^{\I}-$homogeneous, so we will denote them by $x_\beta$, where $\beta\in\N_0^{\I}$ is the degree.

One of the tasks is to check that a relation is primitive in all pre-Nichols algebras with finite $\GK$. The following result will be useful for such proposal.

\begin{remark}\label{rem:coproduct}
Let $I \subset J$ be $\N_0^{\I}$-graded Hopf ideals of $T(V)$, let $S$ be a system of $\N_0^{\I}$-homogeneous generators of $J$, and put  $\toba:=T(V)/I$. Consider an  homogeneous element $\textbf{x}\in J$, and set $\mathbb{Y}_\textbf{x}=\{\textbf{y}\in S:\deg(\textbf{y})<\textbf{x}\}$.
If $\mathbb{Y}_{\textbf{x}} \subset I$, then $\textbf{x}\in\Pc(\toba)$.
\end{remark} 
\pf
Since $J$ is a coideal and the coproduct is $N_0^{\I}$-homogeneous, there exist $a_{\textbf{y}}$, $b_{\textbf{y}}$, $c_{\textbf{y}}$, $d_{\textbf{y}}$, $e_{\textbf{y}}$, $f_{\textbf{y}}\in T(V)$ such that
\begin{align}\label{eq:coproduct-quotients}
\Delta(\textbf{x})=1\ot \textbf{x} + \textbf{x}\ot 1+\sum_{\textbf{y}\in \mathbb{Y}_\textbf{x}}{a_{\textbf{y}}\textbf{y}b_{\textbf{y}}\ot c_{\textbf{y}}+d_{\textbf{y}}\ot e_{\textbf{y}}\textbf{y}f_{\textbf{y}}}\in T(V)\ot T(V).
\end{align}
Then we use that $\textbf{y}=0$ in $\toba$ for all $\textbf{y}\in \mathbb{Y}_{\textbf{x}}$.
\epf

In particular we can apply Remark \ref{rem:coproduct} to $J=\widetilde{\Jc}_{\bq}$, the defining ideal of the distinguished pre-Nichols algebra. This will help us to prove that an element is primitive in a pre-Nichols algebra by studying its image onto $\wtoba_{\bq}$.

Now we provide a refined treatment of some relations already considered in \cite{ASa,ACS}.

\begin{lemma}\label{lem:previous-lemmas} 
Let $\toba$ be a finite $\GK$  pre-Nichols algebra of $\toba_\bq$.
\begin{enumerate}[leftmargin=*,label=\rm{(\alph*)}]
\item \label{item:previous-xij}
Let $i,j\in \I_{\theta}$ be such that $c_{ij}=0$. \emph{Assume  that} one of the following hold: 
\begin{itemize}[leftmargin=*]
\item $\ord q_{ii}+\ord q_{jj}>4$, 
\item $q_{ii}q_{jj}=1$ and there exists $k\in\I_\theta-\{i, j\}$ such that $\qti_{ik}\qti_{jk}\ne 1$. 
\end{itemize}
Then $x_{ij} = 0$ in $\toba$.

\item  \label{item:previous-qsr}
Let $i,j \in \I_{\theta}$ be such that $c_{ij}<0$, $\qti_{ij}^{\,1-c_{ij}}\neq 1$. \emph{Assume  that} one of the following hold: 
\begin{itemize}[leftmargin=*]
\item the Dynkin diagram of $\ku x_i\oplus \ku x_j$ is different from
$\xymatrix @C=15pt{ {\overset{q}{\circ}} \ar  @{-}[r]^{q^{-1}}  & {\overset{q}{\circ} } }, q\in\G_3'$,
\item $c_{ij}=-1$, $q_{ii}=q_{jj}=\qti_{ij}^{-1}\in\G_3'$, and there exists $k\in\I_\theta-\{i,j\}$ such that $\qti_{ik}^2\qti_{jk}\ne 1$.
\end{itemize}
Then $(\adc x_i)^{1-c_{ij}}x_j=0$ in $\toba$.

\item \label{item:previous-qsr-of-distinguished}
Let $i,j\in\I_\theta$ such that $q_{ii}^{1-c_{ij}}=1$. \emph{Assume  that} $\qti_{ij}=q_{ii}$ and one of the following hold 
\begin{itemize}[leftmargin=*]
\item $q_{jj}\ne-1$, 
\item $c_{ij} \le - 2$,
\item $c_{ij} = - 1$, $q_{ii}=-1$ and there exists $k\in \I_\theta- \{i, j\}$ such that $\qti_{jk}, \qti_{ik}^2\qti_{jk}\neq 1$.
\end{itemize}
Then $(\adc x_i)^{1-c_{ij}}x_j=0$ in $\toba$.

\item \label{item:previous-non-Cartan-power}
Let $i\in\I_\theta$ be a non-Cartan vertex. Then $x_i^{N_i}=0$ in $\toba$. 

\item \label{item:previous-xij^2}
Let $i,j\in \I_{\theta}$ be such that $q_{ii}=\qti_{ij}=q_{jj}=-1$ and there exists $k\in\I_{\theta}-\{i,j\}$ such that either $\qti_{ik}^2\neq 1$ or $\qti_{jk}^2\neq 1$. 
\emph{Assume  that} $\qti_{ik}^2\qti_{jk}^2\neq 1$. Then $x_{ij}^2=0$ in $\toba$.

\item \label{item:previous-[x_ijk,x_j]}
Let $i,j,k\in\I_\theta$ be such that $q_{jj}=-1$, $\qti_{ik}=\qti_{ij}\qti_{jk}=1$ and $\qti_{ij} \neq \pm 1$.
\emph{Assume  that} one of the following conditions holds:
\begin{itemize}[leftmargin=*]
\item either $q_{ii}=-1$ or $q_{kk}=-1$,
\item $q_{ii}q_{kk}=1$ and there exists $\ell\in\I_\theta-\{i, j, k\}$ such that $\qti_{i\ell}\neq 1=\qti_{j\ell}=\qti_{k\ell}$,
\item $q_{ii}q_{kk}=1$ and there exists $\ell\in\I_\theta-\{i, j, k\}$ such that $\qti_{j\ell}^2\neq 1=\qti_{i\ell}=\qti_{k\ell}$,
\item $q_{ii}q_{kk}=1$ and there exists $\ell\in\I_\theta-\{i, j, k\}$ such that $\qti_{k\ell}\neq 1=\qti_{j\ell}=\qti_{i\ell}$.
\end{itemize}
Then $[x_{ijk},x_j]_c=0$ in $\toba$.

\item \label{item:previous-[x_{iij},x_{ij}]} Let $i,j\in \I_{\theta}$ be such that $q_{jj}=-1$, $q_{ii}\qti_{ij}\in\G_3'\cup\G_6'$ and either $q_{ii}\in\G_3'$ or $c_{ij} \le - 3$. 
Then $[x_{iij},x_{ij}]_c=0$ in $\toba$.

\item \label{item:previous-[x_{iijk},x_{ij}]} Let $i,j,k\in \I_{\theta}$ be such that $q_{ii}=\pm\qti_{ij}\in\G_3'$, $\qti_{ik}=1$ and either $-q_{jj}=\qti_{ij}\qti_{jk}=1$ or $q_{jj}^{-1}=\qti_{ij}=\qti_{jk}\ne -1$.
Then $[x_{iijk},x_{ij}]_c=0$ in $\toba$.

\item \label{item:previous-triangle-x_{ijk}} Let $i,j,k\in \I_{\theta}$ be such that $\qti_{ij},\qti_{ik},\qti_{jk}\neq 1$. Then 
$$ x_{ijk} - q_{ij}(1-\qti_{jk})x_{j}x_{ik}+\frac{1-\qti_{jk}}{q_{kj}(1-\qti_{ik})}[x_{ik},x_{j}]_c = 0\text{ in }\toba.$$

\item \label{item:previous-[[x_{ij},x_{ijk}],x_j]} Let $i,j,k\in \I_{\theta}$ be such that $q_{ii}=q_{jj}=-1$, $\qti_{ij}^{\,2}=\qti_{jk}^{\,-1}\ne 1$ and $\qti_{ik}=1$. \emph{Assume  that} either 
$q_{kk}^2\ne1$ or $\qti_{ij}^{\, 3}\ne 1$. Then $[[x_{ij},x_{ijk}]_c,x_j]_c=0$ in $\toba$.

\item \label{item:previous-[[x_{ij},[x_{ij},x_{ijk}]],x_j]} Let $i,j,k \in \I_{\theta}$ be such that $q_{ii}=q_{jj}=-1$, $\qti_{ij}^{\ 3}= \qti_{jk}^{\ -1}$ and $\qti_{ik}=1$. 
Then $[[x_{ij},[x_{ij},x_{ijk}]_c]_c,x_j]_c=0$ in $\toba$.

\item \label{item:previous-[[x_{ijk},x_j],x_j]} Let $i,j,k\in\I_{\theta}$ be such that $q_{jj}=\qti_{ij}^{\ 2}=\qti_{jk}\in \G_{3}'$ and $\qti_{ik}=1$. \emph{Assume  that} either $q_{ii}\ne -1$ or $q_{kk}\ne -1$.
Then $[[x_{ijk},x_{j}]_c,x_j]_c=0$ in $\toba$.

\item \label{item:previous-[[[x_{ijk},x_j],x_j],x_j]} 
Let $i,j,k\in \I_{\theta}$ be such that $q_{jj}=\qti_{ij}^{\ 3}=\qti_{jk}\in \G_4'$, $\qti_{ik}=1$. Then $[[[x_{ijk},x_j]_c,x_j]_c,x_j]_c=0$ in $\toba$.

\item \label{item:previous-[x_{ij},x_{ijk}]} 
Let $i,j,k\in \I_{\theta}$ be such that $q_{ii}=\qti_{ij}=-1$, $q_{jj}=\qti_{jk}^{\,-1}\neq -1$ and $\qti_{ik}=1$. Then $[x_{ij},x_{ijk}]_c=0$ in $\toba$.

\item \label{item:previous-[[[x_{ijkl},x_k],x_j],x_k]}
Let $i,j,k,\ell\in \I_{\theta}$ be such that $q_{kk}= -1$, $q_{jj}\qti_{ij}=q_{jj}\qti_{jk}=1$, $\qti_{ik}=\qti_{il}=\qti_{j\ell}=1$ and $\qti_{jk}^2=\qti_{k\ell}^{-1}=q_{\ell\ell}$. Then $[[[x_{ijk\ell},x_{k}]_c,x_{j}]_c,x_{k}]_c=0$ in $\toba$.

\item \label{item:previous-[[[x_{ijk},x_j],[x_{ijkl},x_j]],x_{jk}]} 
Let $i,j,k,\ell \in \I_{\theta}$ and $q\in\Bbbk$ be such that 
$q_{\ell \ell}=\qti_{\ell k}^{\,-1}=q_{kk}=\qti_{jk}^{\,-1}=q^2$ , $\qti_{ij}=q_{ii}^{\,-1}=q^{3}$, $q_{jj}=-1$ and $\qti_{ik}=\qti_{i\ell}=\qti_{j\ell}=1$.
Then $[[[x_{ijk},x_j]_c,[x_{ijk\ell},x_j]_c]_c,x_{jk}]_c=0$ in $\toba$.

\item \label{item:previous-[[x_{ijkl},x_j],x_k]-q_{jk}(q^2-q)[[x_{ijkl},x_k],x_j]} 
Let $i,j,k,\ell \in \I_{\theta}$ be such that $q_{kk}=-1$, $q_{ii}=\qti_{ij}^{\,-1}=q_{jj}^2$, $\qti_{k\ell}=q_{\ell\ell}^{-1}=q_{jj}^{3}$,  $\qti_{jk}=q_{jj}^{-1}$, $\qti_{ik}=\qti_{i\ell}=\qti_{j\ell}=1$.
Then $[[x_{ijk\ell},x_j]_c,x_k]_c-q_{jk}(\qti_{ij}^{\,-1}-q_{jj})[[x_{ijk\ell},x_k]_c,x_j]_c=0$ in $\toba$.

\item \label{item:previous-[x_i,[x_{ijk},x_j]]} 
Let $i,j,k \in \I_{\theta}$ be such that $\ord q_{ii}> 3$, $\qti_{ik}=1$, $q_{jj}=-1$, $\qti_{ij}=q_{ii}^{-2}$, $\qti_{jk}=q_{kk}^{-1}=-q_{ii}^{3}$.
Then $[x_i,[x_{ijk},x_j]_c]_c -\frac{q_{ij}q_{kj}}{1-q_{ii}^{-1}}[x_{ij},x_{ijk}]_c -(q_{ii}+q_{ii}^{2}) q_{ij}q_{ik}x_{ijk}x_{ij}=0$ in $\toba$.
\end{enumerate}
\end{lemma}

\begin{remark}
The extra assumptions imposed on the braiding matrix $\bq$ are in fact necessary. 
If one of these conditions do not hold, then the corresponding relation is not necessarily zero in all finite $\GK$ pre-Nichols algebra, as we see in the following cases:
\begin{itemize}[leftmargin=*]
\item \ref{item:previous-xij}, for $\superqa{3}{q}{\{2\}}$ with $q\in \G_{\infty}$;
\item \ref{item:previous-qsr}, for $A_2$ with $q\in \G_3'$;
\item \ref{item:previous-qsr-of-distinguished}, \ref{item:previous-xij^2} and  \ref{item:previous-[x_ijk,x_j]}, for Cartan $A_{\theta}$ or $D_{\theta}$ with $q=-1$;
\item \ref{item:previous-[x_ijk,x_j]}, for $\superqa{3}{q}{\{1,2,3\}}$, with $q\in \G_{\infty}$;
\item \ref{item:previous-[[x_{ij},x_{ijk}],x_j]}, for $\g(2,3)$ with diagram $d_1$;
\item \ref{item:previous-[[x_{ijk},x_j],x_j]}, for $\g(2,3)$ with diagram $d_2$.
\end{itemize}
\end{remark}

\pf
By \cite[Theorem 1.2]{ACS} if $\bq$ is of Cartan, super or standard type all items are satisfied. The statements \ref{item:previous-xij}, \ref{item:previous-qsr}, \ref{item:previous-qsr-of-distinguished}, \ref{item:previous-non-Cartan-power}, \ref{item:previous-xij^2}, \ref{item:previous-[x_ijk,x_j]}, \ref{item:previous-triangle-x_{ijk}} and \ref{item:previous-[[x_{ijk},x_j],x_j]} were proved in \cite[\S 3]{ACS}.

Items \ref{item:previous-[[x_{ij},[x_{ij},x_{ijk}]],x_j]}, \ref{item:previous-[[[x_{ijk},x_j],x_j],x_j]} and \ref{item:previous-[x_{ij},x_{ijk}]} have corresponding lemmas in \cite[\S 3]{ACS}, but with extra hypothesis on $\bq$. However, we verify by exhaustion that these conditions are satisfied for each diagram in \cite{H-classif}. For 
\ref{item:previous-[x_{iij},x_{ij}]}, \ref{item:previous-[x_{iijk},x_{ij}]}, \ref{item:previous-[[x_{ij},x_{ijk}],x_j]}, \ref{item:previous-[[x_{ij},[x_{ij},x_{ijk}]],x_j]},
\ref{item:previous-[[[x_{ijkl},x_k],x_j],x_k]}, \ref{item:previous-[[[x_{ijk},x_j],[x_{ijkl},x_j]],x_{jk}]}, \ref{item:previous-[[x_{ijkl},x_j],x_k]-q_{jk}(q^2-q)[[x_{ijkl},x_k],x_j]}
and \ref{item:previous-[x_i,[x_{ijk},x_j]]}, it remains to prove that $x_{\beta}$ is primitive. The proof is recursive on the degree $\beta\in\N_0^{\I}$. Indeed we apply Remark \ref{rem:coproduct}
with $J=\widetilde{\Jc}_{\bq}$ and $I$ the ideal generated by those generators of $\widetilde{\Jc}_{\bq}$ of degree $<\beta$.
\epf

\subsection{Verifying more relations}\label{subsec:new-relations}
We study the remaining defining relations of $\wtoba_\bq$ given in  \cite[Theorem 3.1]{An-crelle}.  Here $\toba$ stands for a pre-Nichols algebra of $\bq$ such that $\GK\toba<\infty$.

Let us outline the general strategy.  For each one of the Lemmas below we assume that $\bq$ satisfies the conditions required on \cite[Theorem 3.1]{An-crelle} for including certain relation $x_\beta$ of degree $\beta\in\N_0^{\I}$ in the presentation of $\wtoba_\bq$. Next we suppose that (the image of) $x_\beta$ does not vanish in $\toba$ and prove that (the image of) $x_\beta$ is primitive in $\toba$. Thus we get a braided vector space of diagonal type $V^\bq \oplus \ku x_\beta \subset \Pc(\toba)$ which satisfy $\GK\toba(V^\bq \oplus \ku x_\beta )<\infty$ by Lemma \ref{lem:subspace-primitives}. Now we compute the Dynkin diagram of a suitable chosen subspace of $V^\bq \oplus \ku x_\beta$; this is a straightforward task involving \eqref{eq:notation-qab} and depending only on $\bq$ and $\beta$.
Since we are assuming the validity of Conjecture \ref{conj:AAH}, this diagram should belong to the classification given in \cite{H-classif} and this allow us to arrive at a contradiction. Sometimes we get a Dynkin diagram that belongs to a class in which the conjecture is known to hold true, so we do not need any further assumption.

As in Lemma \ref{lem:previous-lemmas}, the proofs of the lemmas in this subsection use that the relation $x_\beta$ under consideration is primitive in any pre-Nichols algebra $\toba$ with finite $\GK$. This is proved recursively on $\beta$ applying Remark \ref{rem:coproduct} for $J=\widetilde{\Jc}_{\bq}$ and $I$ the ideal generated by those generators of $\widetilde{\Jc}_{\bq}$ of degree $<\beta$.

\begin{lemma}\label{lem:[[x_ij,x_ijk],x_j]}
Let $i,j,k\in \I_{\theta}$ be such that one of the following conditions hold:
\begin{enumerate}[leftmargin=*,label=\rm{(\alph*)}]
\item \label{item:[[x_ij,x_ijk],x_j]-a} $\qti_{ij}=q_{jj}=-1$ and $q_{ii}=-\qti_{jk}^{\,2}\in \G_3'$, $\qti_{ik}=1$;
\item \label{item:[[x_ij,x_ijk],x_j]-b} $q_{kk}=\qti_{jk}=q_{jj}=-1$ and $q_{ii}=-\qti_{ij}\in \G_3'$, $\qti_{ik}=1$; 
\item \label{item:[[x_ij,x_ijk],x_j]-c} $q_{jj}=-1$, $\qti_{ij}=q_{ii}^{-2}\ne 1$, $q_{kk}=\qti_{jk}^{-1}=-q_{ii}^3$, $\qti_{ik}=1$; 
\item \label{item:[[x_ij,x_ijk],x_j]-d} $q_{ii}=q_{jj}=q_{kk}=-1$, $-\qti_{ij}=\qti_{jk}\in\G_3'$, $\qti_{ik}=1$.
\end{enumerate}
Then $[[x_{ij},x_{ijk}]_c,x_j]_c=0$ in $\toba$.
\end{lemma}

\pf
Suppose that $x_\beta:=[[x_{ij},x_{ijk}]_c,x_j]_c\neq 0$ in $\toba$. 

\ref{item:[[x_ij,x_ijk],x_j]-a} We have $\qti_{i\beta}=-q_{ii}\ne 1$, $\qti_{j\beta}=\qti_{jk} \ne 1$,
so the Dynkin diagram of $\ku x_i \oplus \ku x_j\oplus \ku x_\beta$ is a (connected) triangle. Since $\GK\toba<\infty$, Lemma \ref{lem:rank3} implies that $\qti_{ij}\qti_{j\beta}\qti_{i\beta}=1$, but this means $q_{ii}\qti_{jk}=1$ which contradicts $q_{ii}=-\qti_{jk}^{\,2}\in \G_3'$.

\ref{item:[[x_ij,x_ijk],x_j]-b}
In this case $\qti_{i\beta}=-q_{ii}$ and  $\qti_{k\beta}=-1$, so the Dynkin diagram of $\ku x_i \oplus \ku x_j \oplus \ku x_k \oplus \ku x_\beta$ contains a $4$-cycle, which contradicts Lemma \ref{lem:cycles}.

\ref{item:[[x_ij,x_ijk],x_j]-c} Since $\qti_{i\beta}=q_{ii}^{-2}\neq 1$ and  $\qti_{k\beta}=q_{kk}^{-1} \neq 1$, the previous argument applies again.

\ref{item:[[x_ij,x_ijk],x_j]-d} In this case $q_{\beta \beta} = 1$ and $\qti_{i\beta}=-1$, so the Nichols algebra of $\ku x_i \oplus \ku x_\beta$ is infinite GK-dimensional by Lemma \ref{lem:1connected}.
\epf

\begin{lemma}\label{lem:[[xiij,xiijk],xij]} Let $i,j,k\in \I_{\theta}$ be such that $q_{kk}=q_{jj}=\qti_{ij}^{\ -1}=\qti_{jk}^{\ -1}\in \G_9'$, $\qti_{ik}=1$ and $q_{ii}=q_{kk}^6$.
Then $[[x_{iij},x_{iijk}]_c,x_{ij}]_c=0$ in $\toba$.
\end{lemma}
\pf
Assume that $x_\beta :=  [[x_{iij},x_{iijk}]_c,x_{ij}]_c  \ne 0$. 
Since $q_{\beta \beta} = q_{kk}^7$ and $\qti_{k \beta} =q_{kk}^{-1} = (q_{kk}^7)^{-4}$
the braided subspace $W:=\ku x_k \oplus \ku x_\beta \subset \Pc(\toba)$ is of affine Cartan type $A^{(2)}_2$. 
Thus $\GK \toba(W)=\infty$ by \cite[Theorem 1.2 (a)]{AAH-diag}, a contradiction.
\epf

\begin{lemma}\label{lem:[[xijk,xj],xk]} Let $i,j,k\in \I_{\theta}$ be such that $q_{ii}=\qti_{ij}^{\ -1}\in \G_9'$, 
$q_{jj}=\qti_{jk}^{\ -1}=q_{ii}^5$, $\qti_{ik}=1$ and $q_{kk}=q_{ii}^6$. Then $[[x_{ijk},x_{j}]_c,x_k]_c=(1+\qti_{jk})^{-1}q_{jk}[[x_{ijk},x_k]_c,x_j]_c$ in $\toba$.
\end{lemma}
\pf
If $x_\beta =[[x_{ijk},x_{j}]_c,x_k]_c-(1+\qti_{jk})^{-1}q_{jk}[[x_{ijk},x_k]_c,x_j]_c\ne 0$, then $W:=\ku x_k \oplus \ku x_\beta \subset \Pc (\toba)$ has Dynkin diagram 
$\xymatrix@C=30pt{\overset{q}{\underset{k}{\circ}} \ar@{-}[r]^{q} &\overset{q^3}{\underset{\beta}{\circ}}}$ 
where $q :=  q_{ii}^5 \in \G_9'$. As this diagram does not appear in \cite[Table 1]{H-classif}, \cite[Theorem 1.2 (b)]{AAH-diag} assures that $\GK\toba(W)=\infty$, a contradiction with $\GK \toba < \infty$.
\epf

\begin{lemma}\label{lem:[xi,xjjk]} Let $i,j,k\in \I_{\theta}$ be such that $q_{ii}=q_{kk}=-1$, $\qti_{ik}=1$,  $\qti_{ij}\in \G_3'$ and $q_{jj}=-\qti_{jk}=\pm \qti_{ij}$.
Then $[x_{i},x_{jjk}]_c=(1+q_{jj}^2)q_{kj}^{-1}[x_{ijk},x_j]_c+(1+q_{jj}^2)(1+q_{jj})q_{ij}x_jx_{ijk}$ in $\toba$.
\end{lemma}
\pf
Assume that $x_\beta:=[x_{i},x_{jjk}]_c-(1+q_{jj}^2)q_{kj}^{-1}[x_{ijk},x_j]_c-(1+q_{jj}^2)(1+q_{jj})q_{ij}x_jx_{ijk}\neq0$. As $\qti_{\beta i}=\qti_{\beta k}=\qti_{ij}^2$, the Dynkin diagram of $W:=\ku x_i \oplus \ku x_j \oplus  \ku x_k \oplus \ku x_\beta \subset \Pc(\toba)$ contains a $4$-cycle, so $\GK\toba(W)=\infty$ by Lemma \ref{lem:cycles}. This contradicts $\GK \toba <\infty$.
\epf

\begin{lemma}\label{lem:[[xijk,[xjjkl,xk]],xjk]}
Let $i,j,k,\ell\in \I_{\theta}$ be such that 
$\qti_{jk}=\qti_{ij}=q_{jj}^{-1}\in \G_4'\cup \G_6'$, $q_{ii}=q_{kk}=-1$, 
$\qti_{ik}=\qti_{i\ell}=\qti_{j\ell}=1$ 
and $\qti_{jk}^{\ 3}=\qti_{\ell k}$. 
Then $[[x_{ijk},[x_{ijk\ell},x_k]_c]_c,x_{jk}]_c=0$ in $\toba$.
\end{lemma}
\pf
If $x_\beta  :=  [[x_{ijk},[x_{ijk\ell},x_k]_c]_c,x_{jk}]_c \ne 0$, then 
$W=\ku x_i \oplus \ku x_j \oplus \ku x_\beta$ has diagram 
\begin{align*}
\begin{aligned}
\xymatrix@C=30pt@R=10pt{ &\overset{q^3 q_{ll}}{\underset{\beta}{\circ}} & \\ 
\overset{-1}{\underset{i}{\circ}} \ar@{-}[rr]^{q^{-1}} \ar@{-}[ru]^{q^{-1}} & &
\overset{q}{\underset{j}{\circ}}\ar@{-}[lu]_{q^{-5} },} 
\end{aligned}
&& q :=  q_{jj} \in \G_4'\cup \G_6',
\end{align*}
which does not have finite root system by Lemma \ref{lem:rank3}. This contradicts $\GK \toba < \infty$.
\epf

\begin{lemma}\label{lem:[[xijkl,xj],xk]} Let $i,j,k,\ell \in \I_{\theta}$ be such that one of the following hold:
\begin{enumerate}[leftmargin=*,label=\rm{(\roman*)}]
\item \label{item:[[xijkl,xj],xk]-a}
$q_{kk}=-1$, $q_{ii}=\qti_{ij}^{\,-1}=q_{jj}^2$, $\qti_{k\ell}=q_{\ell\ell}^{-1}=q_{jj}^{3}$,  $\qti_{jk}=q_{jj}^{-1}$ and $\qti_{ik}=\qti_{i\ell}=\qti_{j\ell}=1$;
\item \label{item:[[xijkl,xj],xk]-b}
$q_{ii}=\qti_{ij}^{\ -1} = -q_{\ell \ell}^{-1} = -\qti_{kl}$, $q_{jj}=\qti_{jk}=q_{kk}=-1$ and $\qti_{ik}=\qti_{i\ell}=\qti_{j\ell}=1$;
\item \label{item:[[xijkl,xj],xk]-c}
$q_{jj}=\qti_{jk}^{\ -1} \in \G_3'$,  $ q_{ii}=\qti_{ij}^{\ -1}=q_{\ell \ell} = \qti_{kl}^{\ -1}= - q_{j j}$, $q_{kk}=-1$ and $\qti_{ik}=\qti_{i\ell}=\qti_{j\ell}=1$.
\end{enumerate}
Then $[[x_{ijk\ell},x_j]_c,x_k]_c=q_{jk}(\qti_{ij}^{\,-1}-q_{jj})[[x_{ijk\ell},x_k]_c,x_j]_c$ in $\toba$.
\end{lemma}

\pf
Assume that $ x_\beta  :=  [[x_{ijk\ell},x_j]_c,x_k]_c-q_{jk}(\qti_{ij}^{\,-1}-q_{jj})[[x_{ijk\ell},x_k]_c,x_j]_c\ne 0$.

\ref{item:[[xijkl,xj],xk]-a} This is \cite[Lemma 3.24]{ACS}, we included the statement here for completeness.

\ref{item:[[xijkl,xj],xk]-b} Here, $W:=\ku x_i \oplus \ku x_j \oplus \ku x_k\oplus \ku x_\ell\oplus \ku x_\beta \subset \Pc(\toba)$ has Dynkin diagram 
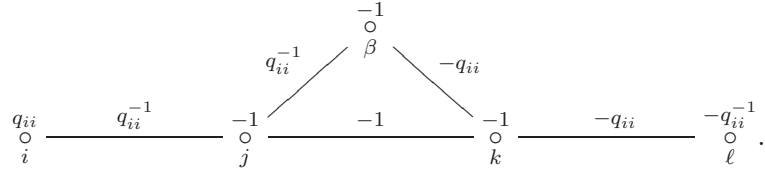
\begin{align*}
\xymatrix@C=30pt@R=10pt{&&&\overset{-1}{\underset{\beta}{\circ}}\ar@{-}[dl]_{q_{ii}^{-1}}\ar@{-}[dr]^{-q_{ii}}&& \\
\overset{q_{ii}}{\underset{i}{\circ}} \ar@{-}[rr]^{q_{ii}^{-1}} && \overset{-1}{\underset{j}{\circ}} \ar@{-}[rr]^{-1}  &&\overset{-1}{\underset{k}{\circ}}\ar@{-}[rr]^{-q_{ii}} && \overset{-q_{ii}^{-1}}{\underset{\ell}{\circ}}. } 
\end{align*}
This diagram does not appear in \cite[Table 4]{H-classif}, a contradiction with $\GK \toba < \infty$.

\ref{item:[[xijkl,xj],xk]-c} Set $W=\ku x_j \oplus \ku x_\beta \subset \Pc (\toba)$, which satisfies $\GK \toba (\ku x_j \oplus \ku x_\beta ) < \infty$ and thus has finite root system \cite[Theorem 1.2 (b)]{AAH-diag}. The Dynkin diagram is
$\xymatrix@C=30pt{\overset{q_{jj}}{\underset{j}{\circ}} \ar@{-}[r]^{-q_{jj}} &\overset{q_{jj}}{\underset{\beta}{\circ}}}$, which does not belong to \cite[Table 1]{H-classif}, a contradiction.
\epf

\begin{lemma}\label{lem:[xi,[xij,xik]]} Let $i,j,k\in \I_{\theta}$ be such that $\qti_{jk}=1$, $q_{ii}=\qti_{ij}=-\qti_{ik}\in \G_3'$. Then $[x_{i},[x_{ij},x_{ik}]_c]_c=-q_{jk}q_{ik}q_{ji}[x_{iik},x_{ij}]_c-q_{ij}x_{ij}x_{iik}$ in $\toba$.
\end{lemma}
\pf
Since $\ku x_j \oplus \ku x_i \oplus \ku x_k$ has finite root system, \cite[Table 2]{H-classif} implies that $q_{jj} = -1$ and $q_{kk} \in \{-1, - q_{ii}^{-1}\}$. 

Assume that $x_\beta  :=  [x_{i},[x_{ij},x_{ik}]_c]_c+q_{jk}q_{ik}q_{ji}[x_{iik},x_{ij}]_c+q_{ij}x_{ij}x_{iik} \ne 0$. The diagram of $\ku x_i\oplus \ku x_\beta$ is 
$d:=\xymatrix@C=30pt{\overset{q_{ii}}{\underset{i}{\circ}} \ar@{-}[r]^{-q_{ii}^2} &\overset{-q_{kk}}{\underset{\beta}{\circ}}}$. 
If $q_{kk}=-1$, then $\GK \toba (\ku x_i\oplus \ku x_\beta) = \infty$ by Lemma \ref{lem:1connected}. In the case $q_{kk}=-q_{ii}^{-1}$, the root system of $d$
is infinite by \cite[Table 1]{H-classif}, hence $\GK \toba (\ku x_i\oplus \ku x_\beta) = \infty$ by \cite[Theorem 1.2 (b)]{AAH-diag}.
\epf

\begin{lemma}\label{lem:[xiijk,xijk]}
Let $i,j,k\in \I_{\theta}$ be such that $q_{jj}=q_{kk}=\qti_{jk}=-1$, $q_{ii}=-\qti_{ij}\in \G_3'$ and $\qti_{ik}=1$. Then $[x_{iijk}, x_{ijk}]_c=0$ in $\toba$.
\end{lemma}

\pf
The degree of $[x_{iijk}, x_{ijk}]_c$ is $\beta :=  3\alpha _i + 2 \alpha _j + 2 \alpha_k$.
Since $q_{\beta \beta}=1$ and $\qti_{i \beta }=q_{ii}^2\ne 1$, it follows from Lemma \ref{lem:1connected} that $[x_{iijk}, x_{ijk}]_c=0$.
\epf

\begin{lemma}\label{lem:[xi,[xij,xj]]}
Let $i,j\in \I_{\theta}$ be such that $-q_{ii},-q_{jj}, \qti_{ij}, q_{ii}\qti_{ij},q_{jj}\qti_{ij}\ne 1$.  Then the relation $[x_i,[x_{ij},x_j]_c]_c=\tfrac{(1+q_{ij})(1-q_{jj}\qti_{ij})}{(1-\qti_{ij})q_{ii}q_{ji}}x_{ij}^2$ holds in $\toba$.
\end{lemma}
\pf
By \cite[Corollary 13]{H-rank3} the diagram $\xymatrix@C=30pt{\overset{q_{ii}}{\underset{i}{\circ}} \ar@{-}[r]^{\qti_{ij}} &\overset{q_{jj}}{\underset{j}{\circ}}}$ 
cannot be extended to a connected diagram of rank $3$ with finite root system. Moreover \cite[Proposition 9 (i)]{H-rank3} warranties that $q_{ii}\qti_{ij}^{\ 2}q_{jj}=-1$ and either $q_{ii}\in\G_3'$ or $q_{jj}\in\G_3'$. By symmetry we can assume $q_{ii}\in \G_3'$. 
If $x_\beta :=[x_i,[x_{ij},x_j]_c]_c-\tfrac{(1+q_{ij})(1-q_{jj}\qti_{ij})}{(1-\qti_{ij})q_{ii}q_{ji}}x_{ij}^2 \ne 0$, then the Nichols algebra of $\ku x_i \oplus \ku x_j \oplus \ku x_\beta \subset \Pc(\toba)$ has finite $\GK$ by Lemma \ref{lem:subspace-primitives}. Hence its diagram must be disconnected by the previous argument. 
But $\qti_{i\beta}= -q_{jj}^{-1}\ne 1$, a contradiction. 
\epf

\begin{remark}\label{rem:rank3-big-roots}
The unique braiding matrix $\bq$ with connected diagram of rank at least $3$, with finite root system and such that $3 \alpha_i+2\alpha_j \in \varDelta_+^\bq$ is
\begin{align}\label{diagram-super-G3}
\xymatrix@C=30pt{\overset{-1}{{\circ}} \ar@{-}[r]^{q^{-1}} & \overset{q}{{\circ}}\ar@{-}[r]^{q^{-3}} &\overset{q^3}{{\circ}},} &&\ord q>3,
\end{align}
which is of type $\mathbf{G}(3)$. Furthermore, this is also the unique $\bq$ such that $m_{ij}\geq 3$ for some $i,j\in\I_3$, see \cite[Table 2]{H-classif}.
Also, for all $\bq$ in \cite[Table 2]{H-classif}, $4 \alpha_i+3\alpha_j$, $5 \alpha_i+3\alpha_j$, $5 \alpha_i+4\alpha_j \notin \varDelta_+^\bq$. We will use these facts frequently in what follows.
\end{remark}

\begin{lemma}\label{lem:diagonal-[xi,x3ai+2aj]}
Assume that the diagram of $\bq$ is one of the following:
\begin{align*}
&\brj(2,3): \xymatrix@C=30pt{\overset{-\zeta^2}{\underset{i}{\circ}}\ar@{-}[r]^{\zeta}&\overset{-1}{\underset{j}{\circ}}}, \, \zeta\in\G_{9}';
&&\ufo(9): \xymatrix@C=30pt{\overset{\zeta}{\underset{i}{\circ}}\ar@{-}[r]^{\zeta^{-5}}&\overset{-1}{\underset{j}{\circ}}}, \, \zeta\in\G_{24}';
\\
&\brj(2,5): \xymatrix@C=30pt{\overset{-\zeta^{-2}}{\underset{i}{\circ}}\ar@{-}[r]^{\zeta^{-2}}&\overset{-1}{\underset{j}{\circ}}}, \, \zeta\in\G_{5}';
&&\text{Standard }\mathbf{G}_2: \xymatrix@C=30pt{\overset{\zeta^{2}}{\underset{i}{\circ}}\ar@{-}[r]^{-\zeta^{-1}}&\overset{-1}{\underset{j}{\circ}}}, \, \zeta\in\G_{8}';
\\
&\ufo(10): \xymatrix@C=30pt{\overset{-\zeta^{-2}}{\underset{i}{\circ}}\ar@{-}[r]^{\pm\zeta^{3}}&\overset{-1}{\underset{j}{\circ}}}, \, \zeta\in\G_{20}';
&&\ufo(11): \xymatrix@C=30pt{\overset{\zeta^{3}}{\underset{i}{\circ}}\ar@{-}[r]^{-\zeta^{4}}&\overset{-\zeta^{-4}}{\underset{j}{\circ}}}, \, \zeta\in\G_{15}';
\\
&\ufo(11): \xymatrix@C=30pt{\overset{\zeta^{3}}{\underset{i}{\circ}}\ar@{-}[r]^{-\zeta^{2}}&\overset{-1}{\underset{j}{\circ}}}, \, \zeta\in\G_{15}';
&&\ufo(12): \xymatrix@C=30pt{\overset{-\zeta^{-2}}{\underset{i}{\circ}}\ar@{-}[r]^{-\zeta^{3}}&\overset{-1}{\underset{j}{\circ}}}, \, \zeta\in\G_{7}'.
\end{align*}
Then $[x_{i},x_{3\alpha_i+2\alpha_j}]_c=\tfrac{1-q_{ii}\qti_{ij}-q_{ii}^2\qti_{ij}^2q_{jj}}{(1-q_{ii}\qti_{ij})q_{ii}}x_{iij}^2$ in $\toba$.
\end{lemma}
\pf
Notice first that either $m_{ij}\in\{4,5\}$ or else $m_{ij}=3$, $q_{jj}=-1$, $q_{ii}\in \G_{4}'$.
For $\bq$ of standard type $\mathbf{G}_2$, the claim was proved in \cite[Lemma 6.8]{ACS}.
For the remaining cases, suppose that $x_{\beta}:=[x_{i},x_{3\alpha_i+2\alpha_j}]_c-\tfrac{1-q_{ii}\qti_{ij}-q_{ii}^2\qti_{ij}^2q_{jj}}{(1-q_{ii}\qti_{ij})q_{ii}}x_{iij}^2\neq 0$.
The diagram of $\Bbbk x_i\oplus \Bbbk x_j\oplus \Bbbk x_{\beta}\subset \Pc(\toba)$ is connected since either $\qti_{j\beta}\neq 1$ when $q_{jj}=-1$, or else $\qti_{i\beta}\neq 1$ when $q_{jj}\ne -1$. This is a contradiction with Remark \ref{rem:rank3-big-roots}.
\epf

\begin{lemma}\label{lem:x4ai+3aj}
Assume that the diagram of $\bq$ is one of the following:
\begin{align*}
&\ufo(7): \xymatrix@C=30pt{\overset{-\overline{\zeta}^2}{\underset{}{\circ}}\ar@{-}[r]^{\overline{\zeta}}&\overset{\zeta}{\underset{}{\circ}}}, \, \zeta\in\G_{12}', &
&\ufo(8): \xymatrix@C=30pt{\overset{-\zeta^2}{\underset{}{\circ}}\ar@{-}[r]^{\zeta^3}&\overset{-1}{\underset{}{\circ}}}, \, \zeta\in\G_{12}',
\\
&\ufo(9): \xymatrix@C=30pt{\overset{\zeta^6}{\underset{}{\circ}}\ar@{-}[r]^{\zeta}&\overset{\overline{\zeta}}{\underset{}{\circ}}}, \, \zeta\in\G_{24}', &
&\ufo(11): \xymatrix@C=30pt{\overset{-\zeta}{\underset{}{\circ}}\ar@{-}[r]^{-\overline{\zeta }^{\,3} }&\overset{\zeta^5}{\underset{}{\circ}} }, \, \zeta \in \G_{15}',
\\
&\brj(2,3): \xymatrix@C=30pt{\overset{-\zeta}{\underset{}{\circ}}\ar@{-}[r]^{\overline{\zeta }^{\,2} }&\overset{\zeta^3}{\underset{}{\circ}} }, \, \zeta \in \G_{9}', &
&\mathbf{G}_2: \xymatrix@C=30pt{\overset{\zeta}{\underset{}{\circ}}\ar@{-}[r]^{-1}&\overset{-1}{\underset{}{\circ}} }, \, \zeta \in \G_{6}'.
\end{align*}
Then $[x_{3\alpha_i+2\alpha_j},x_{ij}]_c=0$ in $\toba$.
\end{lemma}
\pf
Notice first that in all cases $3 \alpha_i+2\alpha_j \in \varDelta_+^\bq$ but $4 \alpha_i+3\alpha_j \notin \varDelta_+^\bq$.
For $\bq$ of type $\mathbf{G}_2$, the claim was proved in \cite[Lemma 4.2]{ACS}. For the five remaining cases, suppose that $x_{\beta}\neq 0$. Consider  $\ku x_i \oplus \ku x_j \oplus  \ku x_{\beta}  \subset \Pc(\toba)$. Since $\tilde{q}_{i\beta}=q_{ii}^8\qti_{ij}^3\ne 1$ one can verify, case-by-case, that the diagram of that subspace is connected. But this contradicts Remark \ref{rem:rank3-big-roots}, since none of these five rank-two diagrams is a subdiagram of \ref{diagram-super-G3}. It must be $x_{\beta}=0$.
\epf

\begin{lemma}\label{lem:diagonal-[xiij,x3ai+2aj]}
Assume that the diagram of $\bq$ is one of the following:
\begin{align*}
&\ufo(11): \xymatrix@C=30pt{\overset{-\zeta^3}{\underset{}{\circ}}\ar@{-}[r]^{-\zeta^4}&\overset{{\overline{\zeta}}^{\,4}}{\underset{}{\circ}}}, \, \zeta \in \G_{15}', &
&\brj(2,3): \xymatrix@C=30pt{\overset{\zeta^3}{\underset{}{\circ}}\ar@{-}[r]^{\overline{\zeta}}&\overset{-1}{\underset{}{\circ}}}, \, \zeta \in \G_{9}'.
\end{align*}
Then $[x_{iij},x_{3\alpha_i+2\alpha_j}]_c=0$ in $\toba$.
\end{lemma}
\pf
In these cases $3 \alpha_i+2\alpha_j \in \varDelta_+^\bq$ but $5 \alpha_i+3\alpha_j \notin \varDelta_+^\bq$.
Supose that $x_{\beta}\neq 0$. If $\bq$ is of type $\brj(2,3)$ then $\qti_{j\beta}=\zeta^4\neq 1$, and for type $\ufo(11)$, $\qti_{i\beta}=-\overline{\zeta}^{\,4}\neq 1$
Hence the diagram of $\ku x_i \oplus \ku x_j \oplus  \ku x_{\beta}  \subset \Pc(\toba)$ is connected, a contradiction with Remark \ref{rem:rank3-big-roots}.
\epf

\begin{lemma}\label{lem:diagonal-[x4ai+3aj,xij]}
Assume that the diagram of $\bq$ is one of the following:
\begin{align*}
&\ufo(10): \xymatrix@C=30pt{\overset{\zeta}{\underset{}{\circ}}\ar@{-}[r]^{\pm\overline{\zeta}^{\,3}}&\overset{-1}{\underset{}{\circ}}}, \, \zeta\in\G_{20}', &
&\ufo(11): \xymatrix@C=30pt{\overset{\zeta^3}{\underset{}{\circ}}\ar@{-}[r]^{-\zeta^4}&\overset{-\overline{\zeta}^{\,4}}{\underset{}{\circ}}}, \, \zeta\in\G_{15}',
\\
&\ufo(11): \xymatrix@C=30pt{\overset{\zeta^5}{\underset{}{\circ}}\ar@{-}[r]^{-\overline{\zeta}^2 }&\overset{-1}{\underset{}{\circ}} }, \, \zeta\in\G_{15}', &
&\brj(2,5): \xymatrix@C=30pt{\overset{\zeta}{\underset{}{\circ}}\ar@{-}[r]^{\zeta^2}&\overset{-1}{\underset{}{\circ}}}, \, \zeta \in \G_{5}'.
\end{align*}
Then $[x_{4\alpha_i+3\alpha_j},x_{ij}]_c=0$ in $\toba$.
\end{lemma}
\pf
Both cases have $4\alpha_i+3\alpha_j \in \varDelta_+^\bq$ and $5\alpha_i+4\alpha_j \not\in \varDelta_+^\bq$. Suppose that $x_{\beta}\neq 0$. 
We check case-by-case that $\qti_{i\beta}=q_{ii}^{10}\qti_{ij}^4\ne 1$, so the diagram of $\ku x_i \oplus \ku x_j \oplus  \ku x_{\beta}  \subset \Pc(\toba)$ is connected, contradicting Remark \ref{rem:rank3-big-roots}. Thus $x_{\beta}=0$.
\epf

\begin{lemma}\label{lem:diagonal-[[xiiij,xiij],xiij]}
Assume that the diagram of $\bq$ is one of the following:
\begin{align*}
&\ufo(10): \xymatrix@C=30pt{\overset{-\overline{\zeta}^2}{\underset{}{\circ}}\ar@{-}[r]^{\pm\zeta^{3}}&\overset{-1}{\underset{}{\circ}}}, \, \zeta\in\G_{20}', &
&\ufo(11): \xymatrix@C=30pt{\overset{\zeta^3}{\underset{}{\circ}}\ar@{-}[r]^{-\zeta^4}&\overset{-\overline{\zeta}^{\,4}}{\underset{}{\circ}}}, \, \zeta \in \G_{15}'.
\end{align*}
Then $[[x_{iiij},x_{iij}]_c,x_{iij}]_c=0$ in $\toba$.
\end{lemma}
\pf
If $x_{\beta}\neq 0$, consider $\ku x_i \oplus \ku x_j \oplus  \ku x_{\beta}  \subset \Pc(\toba)$, which has connected diagram since $\qti_{\beta i}=q_{ii}^{14}\qti_{ij}^3\ne 1$. But this diagram has infinite root system, because it does not appear in \cite[Table 2]{H-classif}, a contradiction.
\epf

\begin{lemma}\label{lem:diagonal-[xiij,x4ai+3aj]}
Assume that the diagram of $\bq$ is one of the following:
\begin{align*}
&\ufo(9): \xymatrix@C=30pt{\overset{-\overline{\zeta}^4}{\underset{}{\circ}}\ar@{-}[r]^{\zeta^{5}}&\overset{-1}{\underset{}{\circ}}}, \, \zeta\in\G_{24}', &
&\ufo(12): \xymatrix@C=30pt{\overset{-\zeta}{\underset{}{\circ}}\ar@{-}[r]^{-\overline{\zeta}^3}&\overset{-1}{\underset{}{\circ}}}, \, \zeta \in \G_{7}'.
\end{align*}
Then $[x_{iij},x_{4\alpha_i+3\alpha_j}]_c=\mathtt{c}_{\bq} \, x_{3\alpha_i+2\alpha_2}^2$ in $\toba$, where $\mathtt{c}_{\bq}\in\Bbbk$ is given in \cite[(3.29)]{An-crelle}.
\end{lemma}
\pf
Supose that $x_{\beta}\neq 0$. As $\qti_{\beta i}=q_{ii}^{12}\qti_{ij}^4\ne 1$, the diagram of $\ku x_i \oplus \ku x_j \oplus  \ku x_{\beta}  \subset \Pc(\toba)$ is connected. Also, $5 \alpha_i+4\alpha_j$ belongs to the set of roots of $\ku x_i \oplus \ku x_j \oplus  \ku x_{\beta}$, so the Nichols algebra of this space has $\GK=\infty$ by Remark \ref{rem:rank3-big-roots}. A contradiction, thus $x_{\beta}= 0$.
\epf

\section{Exceptional cases of type \texorpdfstring{$\g(2,3)$}{}}\label{sec:g23-exceptions}
Theorem \ref{thm:distinguished-is-eminent} says, in particular, that for a braiding matrix $\bq$ of modular, supermodular  or unidentified type, the distinguished pre-Nichols algebra $\wtoba_\bq$  is eminent up to two exceptions. In this section we present, by generators and relations, eminent pre-Nichols algebras $\htoba_{\bq}$ for these two exceptions. We show that in both cases $\htoba_\bq$ fits in an exact sequence of braided Hopf algebras $\Zc \hookrightarrow \htoba_{\bq} \twoheadrightarrow \wtoba_{\bq}$, where $\Zc$ is a $q$-polynomial algebra in two variables. 
Even though the exposition will not make it explicit, the construction of these eminent pre-Nichols algebras was, in some sense, recursive. Namely, we start with a candidate pre-Nichols algebra that covers all finite $\GK$ pre-Nichols. Then we try to show that this candidate has finite $\GK$ by exhibiting a PBW basis; at this point we may realize that some commutation relation is missing. In that case we redefine our candidate, and start again. Luckily, at most two iterations of this process were needed.

The two exceptional diagrams of type  $\g(2,3)$ depend on third-root of unity $\xi$; they are:
\begin{align}\label{diagram:first-exceptional-case-g(2,3)}
\xymatrix@C=30pt{\overset{-1}{\circ} \ar@{-}[r]^{\xi} &\overset{-1}{\circ}\ar@{-}[r]^{\xi}\ &\overset{-1}{\circ}}
\\
\label{diagram:second-exceptional-case-g(2,3)}
\xymatrix@C=30pt{\overset{-1}{\circ} \ar@{-}[r]^{\overline{\xi}} &\overset{\xi}{\circ}\ar@{-}[r]^{\xi}\ &\overset{-1}{\circ}}
\end{align}

\subsection{Type \texorpdfstring{$\g(2,3)$}{}, diagram \texorpdfstring{\eqref{diagram:first-exceptional-case-g(2,3)}}{}}
Fix a braiding matrix $\bq$  with diagram \eqref{diagram:first-exceptional-case-g(2,3)}. The distiguished pre-Nichols algebra has the following presentation:
\begin{align*}
\wtoba_{\bq}=T(V)/\langle x_{1}^2,x_2^2,x_3^2,x_{13}, [[x_{12},x_{123}]_c,x_2]_c,[[x_{123},x_{23}]_c,x_2]_c\rangle
\end{align*}
Notice that Lemma \ref{lem:[[x_ij,x_ijk],x_j]} deals with the last two relations, but under extra assumptions which are not satisfied for this particular $\bq$. We will show, in particular, that there is a pre-Nichols algebra with finite $\GK$ where these elements do not vanish. Set
\begin{align}
x_u&:=[[x_{12},x_{123}]_c,x_2]_c, & x_v&:=[[x_{123},x_{23}]_c,x_2]_c.
\end{align}
As $x_1^2$, $x_2^2$, $x_3^2$, $x_{13}$ are primitive in $T(V)$, they span a Hopf ideal $I:=\langle x_1^2,x_2^2,x_3^2,x_{13}\rangle$ of $T(V)$. Now $x_u,x_v\in \underline{\toba}:= T(V)/I$ are primitive elements by Remark \ref{rem:coproduct}. From \cite[Lemma 2.7]{ASa} we get that $[x_1,x_u]_c$, $[x_u,x_3]_c$, $[x_1,x_v]_c$ and $[x_v,x_3]_c$ are also primitive in $\underline{\toba}$.

\begin{lemma}\label{lem:g(2,3):case-a-[x1,xu]=0}
Let $\toba$ a pre-Nichols algebra of $\toba_{\bq}$ with finite GKdim. Then
\begin{align*}
[x_1,x_u]_c&=[x_v,x_3]_c=[x_u,x_3]_c=[x_1,x_v]_c=0\text{ in }\toba.
\end{align*}
\end{lemma}
\pf
Let $x_{\beta}\in\{ [x_1,x_u]_c, [x_v,x_3]_c, [x_u,x_3]_c, [x_1,x_v]_c\}$. Since $ \underline{\toba} \twoheadrightarrow \toba$ and $x_\beta$ is primitive in $\underline{\toba}$, it is also primitive in $\toba$. Assume $x_{\beta}\ne 0$ in $\toba$. 
By direct computation, $\qti_{1\beta}=\qti_{3\beta}=1$, $\qti_{2\beta}=\xi$ and $q_{\beta\beta}=-1$, so 
the Dynkin diagram of $\ku x_1\oplus \ku x_2\oplus \ku x_3\oplus \ku x_\beta \subset \Pc(\toba)$ is
\begin{align*}
\xymatrix@C=30pt@R=10pt{&&\overset{-1}{\underset{\beta}{\circ}}\ar@{-}[dd]_{\xi}&& \\
&&&&\\
\overset{-1}{\underset{1}{\circ}} \ar@{-}[rr]^{\xi} && \overset{-1}{\underset{2}{\circ}} \ar@{-}[rr]^{\xi}  &&\overset{-1}{\underset{3}{\circ}},} 
\end{align*}
which is not in \cite[Table 3]{H-classif}, contradicting $\GK \toba<\infty$ (we assume Conjecture \ref{conj:AAH}).
\epf

Now we have a candidate for eminent pre-Nichols algebra:
\begin{align}\label{eq:eminent-first-case-g(2,3)}
\htoba_{\bq}=&T(V)/\langle x_{1}^2,x_2^2,x_3^2,x_{13}, [x_1,x_u]_c,[x_1,x_v]_c,[x_u,x_3]_c,[x_v,x_3]_c\rangle.
\end{align}
Notice that this is indeed a braided Hopf algebra, because it is a quotient of the auxiliary $\underline{\toba}$ by an ideal generated by primitive elements.

\begin{prop}\label{prop:g(2,3)-case-a}
Let $\bq$ of type $\g(2,3)$ with Dynkin diagram (\ref{diagram:first-exceptional-case-g(2,3)}). Then
\begin{enumerate}[leftmargin=*,label=\rm{(\alph*)}]
\item \label{item:g(2,3)-case-a-eminent} The pre-Nichols algebra $\htoba_{\bq}$ defined in (\ref{eq:eminent-first-case-g(2,3)}) is eminent, with $\GK \htoba_{\bq} =6$. 
\item  \label{item:g(2,3)-case-a-basis} Consider 
\begin{align*}
x_{12^23^2}&=[x_{123},x_{23}]_c, & 
x_{12^23}&=[x_{123},x_{2}]_c, & 
x_{1^22^33^2}&=[x_{123},x_{12^23}]_c, & 
x_{1^22^23}&=[x_{12},x_{123}]_c.
\end{align*}
Then a basis of $\htoba_{\bq}$ is given by
\begin{align*}
B& =\{x_3^{n_1}x_{23}^{n_2}x_v^{n_3}x_2^{n_4}x_{12^23^2}^{n_5}x_{12^23}^{n_6} x_{1^22^33^2}^{n_7}x_{123}^{n_8}x_u^{n_9}x_{1^22^23}^{n_{10}}x_{12}^{n_{11}}x_1^{n_{12}}:
\\
& \qquad n_1,n_4,n_5,n_7,n_{10},n_{12}\in\{0,1\}, n_i\in \N_0 \text{ otherwise} \}.
\end{align*}
\item  \label{item:g(2,3)-case-a-extension} There is a $\N_0^3$-homogeneous extension of braided Hopf algebras $\Zc \hookrightarrow \htoba_{\bq} \twoheadrightarrow \wtoba_{\bq}$, where  $\Zc$ is the subalgebra of $\htoba_\bq$ generated by $x_u$ and $x_v$. The braided adjoint action of $\htoba_\bq$ on $\Zc$ is trivial, and $\Zc$ is a polynomial algebra in two variables.
\end{enumerate}
\end{prop}

\pf
Lemmas \ref{lem:previous-lemmas} and \ref{lem:g(2,3):case-a-[x1,xu]=0} imply that the projection from $T(V)$ onto each finite $\GK$-dimensional pre-Nichols algebra $\toba$ of $\bq$ factors through $\htoba_{\bq}$. To finish the proof of  \ref{item:g(2,3)-case-a-eminent}, we still need to show that $\GK \htoba_{\bq}=6$. This will be achieved after several steps, where we will simultaneously prove \ref{item:g(2,3)-case-a-basis} and  \ref{item:g(2,3)-case-a-extension}.

\begin{stepo}\label{stepo:xu-xv-no-nulo}
The elements $x_u$ and $x_v$ do not vanish in $\htoba_{\bq}$.
\end{stepo}

\pf
We consider the following representation of $\underline{\toba}$, $\rho:\underline{\toba}\longrightarrow \Bbbk^{4\times 4}$,
\begin{align*}
\rho(x_1)& = \left(\begin{smallmatrix}
0 & 0 & 1 & 0 \\ 0 & 0 & 0 & 1 \\ 0 & 0 & 0 & 0 \\ 0 & 0 & 0 & 0
\end{smallmatrix}\right),
&
\rho(x_2)& = \left(\begin{smallmatrix}
0 & 0 & 0 & 1 \\ 0 & 0 & 0 & 0 \\ 0 & 0 & 0 & 0 \\ 0 & 0 & 0 & 0
\end{smallmatrix}\right),
&
\rho(x_3)&=\left(\begin{smallmatrix}
0 & 1 & 0 & 0 \\ 0 & 0 & 0 & 0 \\ 0 & 0 & 0 & q_{13} \\ 0 & 0 & 0 & 0
\end{smallmatrix}\right).
\end{align*}
Then $\rho(x_u)\ne 0$, $\rho(x_v)\ne 0$ since the first rows of $\rho(x_u),\rho(x_v)$ are not zero; hence $x_u,x_v\ne 0$ in $\underline{\toba}$. As $\underline{\toba}^n=\htoba_{\bq}^n$ if $n\leq 6$, we have that $x_u,x_v\ne 0$ in $\htoba_{\bq}$.
\epf

\begin{stepo}
The adjoint action of $\htoba_\bq$ on $\Zc$ is trivial, and $\Zc$ has basis $\{x_u^m x_v^n: m,n\in\N_0\}$.
\end{stepo}
\pf 
By Step \ref{stepo:xu-xv-no-nulo}, $x_u,x_v\ne0$. Moreover $x_u^n,x_v^n\ne 0$ for all $n\in\N$ since they are primitive elements such that $q_{uu}=q_{vv}=1$, and
$\Zc$ is a $q$-polynomial algebra in the variables $x_u, x_v$.

By definition of $\htoba_{\bq}$, $(\ad_c x_i)x_u=(\ad_c x_i)x_v=0$ for $i=1,3$, and $[x_u,x_2]_c=[x_v,x_2]_c=0$ since $x_2^2=0$ in $\htoba$. So $(\ad_c x)x_u=(\ad_c x)x_v=0$ for every homogeneous element $x\in\htoba_{\bq}$ of positive degree.
\epf

\begin{stepo}
The linear span of $B$ is $\htoba_{\bq}$.
\end{stepo}
\pf
It is enough to check that the subspace $I$ spanned by $B$ is a left ideal of $\htoba_{\bq}$. 
From $[x_1,x_u]_c=[x_v,x_3]_c=0$ and \eqref{eq:braided-commutator-iteration} we get the equalities:
\begin{align*}
x_{12}^3 x_3 &= q_{13}^3q_{23}^3 x_3 x_{12}^3, &
x_{1} x_{23}^3 &= q_{12}^3q_{13}^3 x_{23}^3 x_1.
\end{align*}
From these equalities we obtain the following:
\begin{align}\label{eq:g(2,3)-a-relations1}
[x_{12}, x_{1^22^23}]_c&=0, & [x_{12^23^2},x_{23}]_c &=0.
\end{align}
Using \eqref{eq:braided-commutator-iteration} again and $x_{13}=x_1^2=x_3^2=0$ we also get
\begin{align*}
[x_{23},x_3]_c&=[x_{123},x_3]_c=0, & [x_{1^22^23}, x_3]_c&=\zeta^2q_{13}q_{23} x_{123}^2,
\\
[x_1,x_{12}]_c&=[x_1,x_{123}]_c=0, & [x_1,x_{12^23^2}]_c&=\zeta^2q_{12}q_{13} x_{123}^2.
\end{align*}
Using the last equality, together with $[x_u,x_3]_c=0$, \eqref{eq:braided-commutator-right-mult} and \eqref{eq:braided-commutator-iteration}, we get
\begin{align*}
0 &= \left[ x_1, [x_u,x_3]_c \right]_c=\left[ x_1, [[x_{1^22^23},x_2]_c,x_3]_c \right]_c
\\
&= \left[ x_1, [x_{1^22^23},x_{23}]_c +q_{12}^2q_{13}x_2x_{123}^2- \zeta^2q_{13}q_{23}^2 x_{123}^2x_2\right]_c
\\
&= q_{12}^2q_{13} x_{1^22^23}x_{123}+\zeta q_{12}^2q_{13}^2q_{23}x_{123}x_{1^22^23} +q_{12}^2q_{13}x_{12}x_{123}^2- \zeta^2q_{12}^2q_{13}^3q_{23}^2 x_{123}^2x_{12}
\\
&= 2q_{12}^2q_{13}x_{12}x_{123}^2-2\zeta^2 q_{13}^3q_{23}^2q_{12}^2 x_{123}^2x_{12};
\end{align*}
thus $x_{12}x_{123}^2=\zeta^2q_{13}^2q_{23}^2 x_{123}^2x_{12}$. Analogously, $x_{123}^2x_{23}=\zeta^2q_{12}^2q_{13}^2 x_{23}x_{123}^2$. From these two equalities we deduce the following:
\begin{align*}
[x_{1^22^23},x_{123}]_c&=0, & [x_{123},x_{12^23^2}]_c&=0.
\end{align*}
Using the first equation and the first one of \eqref{eq:g(2,3)-a-relations1}.
\begin{align*}
x_{1^22^23}^2=x_{1^22^23}(x_{12}x_{123}-\zeta q_{13}q_{23}x_{123}x_{12})=-(x_{12}x_{123}-\zeta q_{13}q_{23}x_{123}x_{12})x_{1^22^23}=-x_{1^22^23}^2.
\end{align*}
This computation and an analogous one for $x_{12^23^2}$ imply that $x_{1^22^23}^2=x_{12^23^2}^2=0$.

Next we check the following equations:
\begin{align*}
[x_1, x_{12^23}]_c &= q_{12}q_{32} x_{1^22^23}+(\zeta-1)q_{13}q_{23}x_{123}x_{12},
\\
[x_1, x_{1^22^33^2}]_c &= -q_{12}q_{13}x_{123}[x_1, x_{12^23}]_c+\zeta^2q_{12}q_{32}[x_1, x_{12^23}]_c x_{123}=0.
\end{align*}
In a similar way we get the equalities:
\begin{align*}
[x_{12^23},x_2]_c &= 0, & [x_{1^22^33^2},x_2]_c &=(1-\zeta^2)q_{12}q_{32}x_{12^23}^2,
\\
[x_{12^23},x_3]_c &=x_{1^22^23}, & [x_{1^22^33^2},x_3]_c &=0.
\end{align*}
Using the relations involving $x_{1^22^33^2}$ we obtain that $x_{1^22^33^2}^2=0$.

A rutinary recursive proof shows that $x_{\alpha}x_{\beta}=q_{\alpha\beta} x_{\beta}x_{\alpha}+$ ordered products of intermediate PBW generators for each pair of roots $\alpha<\beta$, so the step is proved.
\epf

\begin{stepo}
There is a degree-preserving extension of braided Hopf algebras $\Zc \hookrightarrow \htoba_{\bq} \twoheadrightarrow \wtoba_{\bq}$. Furthermore $B$ is a basis of $\htoba_{\bq}$ and $\GK \htoba_{\bq}=\GK\wtoba_{\bq}+\GK \Zc=6$.
\end{stepo}
Let $\Zc'=\htoba^{\co \pi}$ for $\htoba_{\bq} \twoheadrightarrow \wtoba_{\bq}$. Since $\Zc\subseteq \Zc'$, from \cite[Lemma 2.4]{ACS} we get
\begin{align*}
\mathcal{H}_{\htoba_{\bq}}&=\mathcal{H}_{\Zc'}\mathcal{H}_{\wtoba_{\bq}}\geq \mathcal{H}_{\Zc}\mathcal{H}_{\wtoba_{\bq}}\geq\\
&\geq \frac{1}{(1-t_1^2t_2^3t_3)(1-t_1t_2^3t_3^2)}\frac{(1+t_1)(1+t_2)(1+t_1^2t_2^2t_3)(1+t_1^2t_2^3t_3^2)(1+t_1t_2^2t_3^2)(1+t_3)}{(1-t_1t_2)(1-t_1t_2^2t_3)(1-t_2t_3)(1-t_2t_3)}.
\end{align*}
On the other hand $\htoba_{\bq}$ is spanned by $B$, so
\begin{align*}
\mathcal{H}_{\htoba_{\bq}}\leq \frac{(1+t_1)(1+t_2)(1+t_1^2t_2^2t_3)(1+t_1^2t_2^3t_3^2)(1+t_1t_2^2t_3^2)(1+t_3)}{(1-t_1t_2)(1-t_1t_2^2t_3)(1-t_2t_3)(1-t_2t_3)(1-t_1^2t_2^3t_3)(1-t_1t_2^3t_3^2)}
\end{align*}
These inequalities between the Hilbert series say that
\begin{align*}
\mathcal{H}_{\htoba_{\bq}}= \frac{(1+t_1)(1+t_2)(1+t_1^2t_2^2t_3)(1+t_1^2t_2^3t_3^2)(1+t_1t_2^2t_3^2)(1+t_3)}{(1-t_1t_2)(1-t_1t_2^2t_3)(1-t_2t_3)(1-t_2t_3)(1-t_1^2t_2^3t_3)(1-t_1t_2^3t_3^2)}
\end{align*}
so $\Zc=\Zc'$, $B$ is a basis of $\htoba_{\bq}$ and $\GK \htoba_{\bq}=\GK\wtoba_{\bq}+\GK \Zc=4+2=6$.
\epf

\subsection{Type \texorpdfstring{$\g(2,3)$}{}, diagram \texorpdfstring{\eqref{diagram:second-exceptional-case-g(2,3)}}{}}

Let $\bq$ be a braiding matrix with Dynkin diagram \eqref{diagram:second-exceptional-case-g(2,3)}.
In this case the distiguished pre-Nichols algebra is
\begin{align*}
\wtoba_{\bq}=T(V)/\langle x_{1}^2,x_3^2,x_{13}, [x_{223},x_{23}]_c,x_{221},x_{2223}, [[x_{123},x_2]_c,x_2]_c \rangle.
\end{align*}
Notice that Lemma \ref{lem:previous-lemmas} \ref{item:previous-[[x_{ijk},x_j],x_j]} deals with the relation $x_u:=[[x_{123},x_2]_c,x_2]_c$ but under extra assumptions which are not satisfied for this particular $\bq$. We will see that $x_u$ is not zero in at least one pre-Nichols algebra with finite $\GK$. Set:
\begin{align}\label{eq:eminent-second-case-g(2,3)-PBWgen}
x_{12^23}&=[x_{123},x_{2}]_c, \, \,
x_{12^23^2}=[x_{123},x_{23}]_c, \,    \,
x_{12^33^2}=[x_{12^23^2},x_2]_c, \,  \,
x_v=[x_{123},x_{12^23}]_c.
\end{align}
Consider the algebra:
\begin{align}\label{eq:eminent-second-case-g(2,3)}
\htoba_{\bq}=T(V)/\langle x_{1}^2,x_3^2,x_{13}, [x_{223},x_{23}]_c,x_{221},x_{2223},[x_{v},x_3]_c,[x_{12^33^2},x_2]_c,[x_{12^33^2},x_3]_c\rangle.
\end{align}
Next we prove that $\htoba_{\bq}$ is an eminent pre-Nichols algebra.

\begin{prop}\label{prop:g(2,3)-case-b}
Let $\bq$ is of type $\g(2,3)$ with Dynkin diagram (\ref{diagram:second-exceptional-case-g(2,3)}). Then
\begin{enumerate}[leftmargin=*,label=\rm{(\alph*)}]
\item \label{item:g(2,3)-case-b-eminent} The algebra $\htoba_{\bq}$ defined in \eqref{eq:eminent-second-case-g(2,3)} is an eminent pre-Nichols of $\bq$, with $\GK \htoba_{\bq} =6$. 
\item  \label{item:g(2,3)-case-b-basis} A basis of $\htoba_{\bq}$ is given by
\begin{align*}
B=\{x_3^{n_1}x_{23}^{n_2}x_{223}^{n_3}x_2^{n_4}x_{12^33^2}^{n_5}x_{12^23^2}^{n_6}x_u^{n_7}x_{12^23}^{n_8}x_{123}^{n_9}
x_v^{n_{10}}
x_{12}^{n_{11}}x_1^{n_{12}}:n_1,n_3,n_5,n_6,n_{10},n_{11}\in\{0,1\}\}
\end{align*}
\item  \label{item:g(2,3)-case-b-extension} There is a $\N_0^3$-homogeneous extension of braided Hopf algebras $\Zc \hookrightarrow \htoba_{\bq} \twoheadrightarrow \wtoba_{\bq}$, where  $\Zc$ is the subalgebra of $\htoba_\bq$ generated by $x_u$ and $x_v$. The braided adjoint action of $\htoba_\bq$ on $\Zc$ is trivial, and $\Zc$ is a polynomial algebra in two variables.
\end{enumerate}
\end{prop}
\pf
We proceed in several steps. The first two steps are devoted to verify that the defining ideal of $\htoba_\bq$ is a Hopf ideal, and also that $\htoba_\bq$ project onto an arbitrary pre-Nichols algebra $\toba$ with finite $\GK$. Consider the following auxiliary algebra:
$$ \underline{\toba}:= T(V)/\langle x_1^2, x_3^2, x_{13}, x_{221}, x_{2223}, [x_{223},x_{23}]_c,[x_{v},x_3]_c\rangle.$$

\begin{stepi}\label{stepi:Bbarra-Hopf}
$\underline{\toba}$ is a braided Hopf algebra and the canonical projection $T(V)\rightarrow \toba$ induces a surjective Hopf algebra map $\pi:\underline{\toba}\rightarrow \toba$.
Also, $x_u$ and $x_v$ are primitive.
\end{stepi}
\pf
Let $\toba'=T(V)/\langle x_1^2, x_3^2, x_{13}, x_{221}, x_{2223}, [x_{223},x_{23}]_c\rangle$.
As $x_1^2$, $x_3^2$, $x_{13}$, $x_{221}$, $x_{2223}$ are primitive in $T(V)$, $J=\langle x_1^2, x_3^2, x_{13}, x_{221}, x_{2223}\rangle$ is a Hopf ideal of $T(V)$. Also, $[x_{223},x_{23}]_c\in T(V)/J$ is primitive by Remark \ref{rem:coproduct}, so $\toba'$ is a braided Hopf algebra.
By Lemma \ref{lem:previous-lemmas}, the canonical projection $T(V)\twoheadrightarrow \toba$ induces a surjective Hopf algebra map $\toba'\twoheadrightarrow \toba$.

Next we prove that $x_v$ is primitive in $\toba'$. 
By \eqref{eq:coproduct-quotients} applied to $\toba'\twoheadrightarrow \wtoba_{\bq}$,
\begin{align*}
\Delta(x_v) \in 1\ot x_v+x_v\ot 1 + \toba'\ot \langle x_u\rangle+\langle x_u\rangle\ot \toba'.
\end{align*}
Also, $x_u$ is primitive in $\toba'$ by Remark \ref{rem:coproduct}, and using \texttt{GAP},
\begin{align}\label{eq:g(2,3)-b-xu-x1-x3}
[x_u,x_1]_c &=[x_u,x_3]_c =0.
\end{align}
Notice that $x_u$ and $x_v$ are the superletters associated to the Lyndon words $x_1x_2x_3x_2^2$ and $x_1x_2x_3x_1x_2x_3x_2$, respectively, according to the definitions in \cite{Kh}.
By \cite[Lemma 13]{Kh} and these relations, there exist $a,b,c\in \Bbbk$ such that
$$\Delta(x_v)=1\ot x_v+x_v\ot 1+ax_1\ot x_ux_3+bx_1x_u\ot x_3+cx_1x_3\ot x_u.$$
As $(\Delta\ot \id)\Delta(x_v)=(\id\ot \Delta)\Delta(x_v)$ we have that $a=b=c=0$. 
Now, $[x_{v},x_3]_c$ is primitive in $\toba'$ by \cite[Lemma 2.7]{ASa} so $\underline{\toba}$ is a braided Hopf algebra.

Finally, suppose that $x_{\beta}:=[x_{v},x_3]_c\ne 0$ in $\toba$. 
By inspection, if the diagram of a matrix $\bq'$ is connected and contains \eqref{diagram:second-exceptional-case-g(2,3)},
then $\bq'$ is not in \cite[Table 3]{H-classif}; thus $\GK\toba_{\bq'}=\infty$, assuming Conjecture \ref{conj:AAH}.
Now, the diagram of $\ku x_1\oplus \ku x_2\oplus \ku x_3\oplus \ku x_\beta$ is connected since $\qti_{2\beta}=\zeta$, and 
we get a contradiction.
\epf

\begin{stepi}\label{stepi:g(2,3):case-b-new-relations}
$[x_{12^33^2},x_2]_c$, $[x_{12^33^2},x_3]_c$ are primitive in $\underline{\toba}$, thus $\htoba_{\bq}$ is a braided Hopf algebra.
The canonical projection $T(V)\rightarrow \toba$ induces a surjective Hopf algebra map $\pi:\htoba_{\bq}\rightarrow \toba$.
\end{stepi}
\pf
From $x_{221}=x_{2223}=0$ we get $x_1x_{2}^3=q_{12}^3x_2^3x_1$ and $x_3x_{2}^3=q_{32}^3x_2^3x_3$. From the last two equalities we deduce the following:
\begin{align}\label{eq:g(2,3)-b-xu-x2}
[x_u,x_2]_c&=[[[x_{123},x_2]_c,x_2]_c,x_2]_c=x_{123}x_2^3-q_{12}^3q_{32}^3x_2^3x_{123}=0.
\end{align}
As in the proof of Step \ref{stepi:Bbarra-Hopf}, we use \cite[Lemma 13]{Kh} to show that $[x_{12^33^2},x_2]_c$ and $[x_{12^33^2},x_3]_c$ are primitive,
hence $\htoba_{\bq}$ is a braided Hopf algebra since $\htoba_{\bq}=\underline{\toba}/ \langle [x_{12^33^2},x_2]_c,[x_{12^33^2},x_3]_c \rangle$.

Assume that $x_\beta =[x_{12^33^2},x_i]_c \ne 0$, $i\in\I_{2,3}$. The diagram of $\ku x_1\oplus \ku x_2\oplus \ku x_3\oplus \ku x_\beta$ is connected since $\qti_{2\beta}\ne 1$; 
the same argument as in Step \ref{stepi:Bbarra-Hopf} leads to a contradiction  since we assume Conjecture \ref{conj:AAH}.
\epf

By Step \ref{stepi:g(2,3):case-b-new-relations}, it is enough to prove that $\GK \htoba_{\bq} <\infty$. 
To do so, we will see that $B$ is a basis of $\htoba_{\bq}$ in three steps.

\begin{stepi}\label{stepi:Z-normal}
The adjoint action of $\htoba_\bq$ on $\Zc$ is trivial, and $\Zc$ has basis $\{x_u^m x_v^n: m,n\in\N_0\}$.
\end{stepi}

\pf
Using \texttt{GAP} we check that $x_u, x_v\ne 0$ in the pre-Nichols algebra $\toba'$ introduced in the proof of Step \ref{stepi:Bbarra-Hopf}; 
hence $x_u, x_v\ne 0$ in $\htoba_{\bq}$ since $(\toba')^{\alpha}=\htoba_{\bq}^{\alpha}$ for all $\alpha \le 2\alpha_1+3\alpha_2+2\alpha_3$. 
Therefore, $x_u^n, x_v^n\ne 0$ for all $n\in\N$ since $x_u$ and $x_v$ are primitive and $q_{uu}=q_{vv}=1$; moreover $\{x_u^mx_v^n: m,n\in\N_0\}$ is a basis of $\Zc$, and $\Zc$ is a Hopf subalgebra. By \eqref{eq:g(2,3)-b-xu-x1-x3} and \eqref{eq:g(2,3)-b-xu-x2}, $(\ad_c x_i)x_u=0$ for all $i\in\I_3$; hence $(\ad_c x)x_u=0$ for all $x\in\htoba_{\bq}$ homogeneous of degree $>0$. 
\epf

\begin {stepi}
$\htoba_{\bq}$ is spanned by $B$.
\end {stepi}
\pf
To prove the statement, we will see that the subspace $I$ spanned by $B$ is a left ideal of $\htoba_{\bq}$.
As $x_1^2=x_3^2=0$, we also have: 
\begin{align*}
x_1x_{12}&=-q_{12}x_{12}x_1,& x_1x_{123}&=-q_{12}q_{13}x_{123}x_1,\\
x_{23}x_3&=-q_{23}x_3x_{23},& x_{123}x_3&=-q_{13}q_{23}x_3x_{123}.
\end{align*}
From $(\ad_c x_2)^3x_3=[x_{223},x_{23}]_c=0$ we deduce the following equality:
\begin{align*}
x_{223}^2&=x_{223}(x_2x_{23}-\xi q_{23}x_{23}x_2)=-\xi^{-2}(x_2x_{23}-\xi q_{23}x_{23}x_2)x_{223}=-\xi^{-2}x_{223}^2.
\end{align*}
Hence $x_{223}^2=0$.

From $x_{221}=0$ we also have $[x_{12},x_2]_c=0$. Using this relation and $x_1^2=0$ we check that
$x_{12}^2=0$; therefore,
\begin{align*}
x_{12}x_{123}= x_{12}(x_{12}x_3-q_{13}q_{23}x_3x_{12})=-q_{13}q_{23}(x_{12}x_3-q_{13}q_{23}x_3x_{12})x_{12}=-q_{13}q_{23}x_{123}x_{12}.
\end{align*}
Using the relations already proved, \eqref{eq:iterated-adjoint}, \eqref{eq:eminent-second-case-g(2,3)-PBWgen} and (\ref{eq:braided-commutator-iteration}), the following relations also hold: 
\begin{align*}
[x_{12},x_{12^23}]_c&=0, & [x_{12^23^2},x_3]_c&=0, & [x_{12^33^2},x_{23}]_c&=0,
\\
[x_{12^33^2},x_{223}]_c&=0, & [x_{12^23^2},x_{23}]_c&=0, & [x_{12^23^2},x_{223}]_c&=0.
\end{align*}
Similarly,
\begin{align*}
[x_{223},x_3]_c&=\xi^2q_{23}x_{23}^2, & [x_{12^23},x_3]_c&=x_{12^23^2},
\\
[x_1,x_{12^23}]_c&=(\xi^2-1) q_{12}q_{13}x_{123}x_{12}, & [x_1,x_{12^23^2}]_c&=(\xi^2-1)q_{12}q_{13}x_{123}^2,
\\
[x_{12},x_{23}]_c&=(\xi-1)q_{12}x_2x_{123}-\xi q_{23}x_{12^23}, &[x_{12^23},x_{23}]_c&=-q_{23}x_{12^33^2}-q_{12}q_{32}x_2x_{12^23^2}.
\end{align*}
Next we use \eqref{eq:braided-commutator-iteration} and $[x_v,x_3]_c=0$ to deduce that $[x_{123},x_{12^23^2}]_c=0$; also,
\begin{align*}
x_{12^23^2}^2&=(x_{123}x_{23}-q_{123,23}x_{23}x_{123})x_{12^23^2}
=-x_{12^23^2}(x_{123}x_{23}-q_{123,23}x_{23}x_{123})=-x_{12^23^2}^2.
\end{align*}
Hence $x_{12^23^2}^2=0$. From here we check that $[x_{12^23^2},x_{12^33^2}]_c=0$; this relation and $[x_{12^33^2},x_2]_c=0$ imply that $x_{12^33^2}^2=0$.

Again, we prove recursively that $x_{\alpha}x_{\beta}=q_{\alpha\beta} x_{\beta}x_{\alpha}+$ ordered products of intermediate PBW generators for each pair of roots $\alpha<\beta$, so the step is proved.
\epf

\begin{stepi}
There is a degree-preserving extension of braided Hopf algebras $\Zc \hookrightarrow \htoba_{\bq} \twoheadrightarrow \wtoba_{\bq}$. 
Furthermore $B$ is a basis of $\htoba_{\bq}$ and $\GK \htoba_{\bq}=\GK\wtoba_{\bq}+\GK \Zc=6$.
\end{stepi}
The proof is analogous to the corresponding step in Proposition \ref{prop:g(2,3)-case-a}.
Indeed, Step \ref{stepi:Z-normal} shows that $\Zc$ is a central Hopf subalgebra of $\htoba_{\bq}$ with basis $\{x_u^mx_v^n: m,n\in\N_0\}$.
If $\Zc':=\htoba^{\co \pi}$, then $\Zc\subseteq \Zc'$ and \cite[Lemma 2.4]{ACS} implies that
$\mathcal{H}_{\htoba_{\bq}}=\mathcal{H}_{\Zc'}\mathcal{H}_{\wtoba_{\bq}}\geq \mathcal{H}_{\Zc}\mathcal{H}_{\wtoba_{\bq}}$. 
On the other hand $\htoba_{\bq}$ is spanned by $B$, so we have an equality between the Hilbert series:
\begin{align*}
\mathcal{H}_{\htoba_{\bq}}=\mathcal{H}_{\Zc}\mathcal{H}_{\wtoba_{\bq}} 
=\frac{(1+t_1)(1+t_1t_2)(1+t_1t_2^2t_3^2)(1+t_1t_2^3t_3^2)(1+t_2^2t_3)(1+t_3)}
{(1-t_1t_2t_3)(1-t_1t_2^2t_3)(1-t_1^2t_2^3t_3^2)(1-t_1t_2^3t_3)(1-t_2)(1-t_2t_3)}.
\end{align*}
Thus $\Zc=\Zc'$, $B$ is a basis of $\htoba_{\bq}$ and $\GK \htoba_{\bq}=\GK\wtoba_{\bq}+\GK \Zc=6$.
\epf

\end{document}